\documentclass[11pt, reqno]{amsart}

\usepackage{amsfonts}
\usepackage{amsmath}
\usepackage[latin1]{inputenc}
\usepackage[english]{babel}
\usepackage{amsthm}
\usepackage{graphicx}
\usepackage{subfigure}
\usepackage{color}

\usepackage{amssymb}
\usepackage{amsopn}
\usepackage{enumerate}
\usepackage{comment}
\usepackage{wrapfig}

\usepackage{pifont,bbding}
\usepackage{comment}
\usepackage{mathrsfs}

\newcommand{\N}{\mathbb{N}}
\newcommand{\R}{\mathbb{R}}

\newcommand{\A}{\mathcal{A}}

\renewcommand{\P}{P_\alpha}

\renewcommand{\S}{\mathscr S}

\newcommand{\tr}{\mathrm{tr}}
\newcommand{\E}{E_{\mathrm{isop}}}  

\newcommand{\mP}{\mathscr P_\alpha}

\theoremstyle{plain}
\newtheorem{thm}{Theorem}[section]
\newtheorem{prop}[thm]{Proposition}

\newtheorem{lem}[thm]{Lemma}
\theoremstyle{definition}
\newtheorem{defn}[thm]{Definition}

\theoremstyle{remark}
\newtheorem{rem}[thm]{Remark}

\numberwithin{equation}{section}

\title[Minimal partition with trace constraint]
{A minimal partition problem with trace constraint\\ in the Grushin plane.}

\author[V.~Franceschi]{Valentina Franceschi}
\email{valentina.franceschi@inria.fr}

\address
{Inria, team GECO \& CMAP, \'Ecole Polytechnique, CNRS, Universit\'e
Paris-Saclay, Palaiseau, France}


\begin{document}

\begin{abstract}
We study a variational problem for the perimeter associated with the Grushin plane, called minimal partition problem with trace constraint.
This consists in studying how to enclose three prescribed
areas in the Grushin plane, using the least amount of perimeter, under an additional 
``one-dimensional" constraint on the intersections of their boundaries.
We prove existence of regular solutions for this problem, and we characterize them in terms of isoperimetric sets, showing differences with the Euclidean case.
The problem arises from the study of quantitative isoperimetric inequalities and has connections with the theory of minimal clusters. 
\end{abstract}

\maketitle

\section{Introduction}

Minimal partition problems are variational problems in which the aim 
is to find the best configuration of a given number of regions in the 
space in order to minimize their total perimeter.
For a volume measure $V$ and a perimeter measure $\mathcal P$ on $\R^n$, $n\in\N$,
a {\em cluster}
is a finite disjoint family of sets $\{E_i\}_{i=1}^m$, $E_i\subset\R^n$ 
having finite perimeter $\mathcal P$ and finite, positive volume $V$.  
Given $v_i>0$ for $i=1,\dots,m$, by a {\em minimal partition problem} we mean any problem of the form
\[
\inf\Big\{\mathscr P_{\mathcal P}(E) : E=E_1\cup\dots\cup E_m\subset\R^2,\  V(E_i)=v_i\Big\}
\]
where $\{E_i\}_{i=1}^m$ is a cluster and
$\mathscr P_\mathcal{P}$ is defined by 
\begin{equation}
\label{eq:mP}
\mathscr P_{\mathcal P}(E)=\frac{1}{2}\Big\{\mathcal P(\R^n\setminus E)+\sum_{i=1}^m\mathcal P(E_i)\Big\}
\end{equation}
(cf. \cite{MagGMT, Morg}). When $\mathcal P=P$ is the De Giorgi perimeter and $V=\mathcal L^n$ is the $n$-dimensional Lebesgue measure, existence of minimizers is proved in \cite{Alm}.
In the 1993 paper \cite{Foisy93}, the authors provide the first complete solution to a minimal partition problem, in the case when $n=2$, $\mathcal P=P$, $V=\mathcal L^n$ and $m=2$: 
the unique minimizers are the so called {\em double bubbles}. 
This result has been extended to $n\geq 3$ in \cite{HMRR02, Reic},
while for 
general dimensions
 and $m\geq 3$, several open questions about minimal clusters in the Euclidean setting are still open (see \cite{Morg}).

In the present paper, we study a planar minimal partition problem {\em with trace constraint} that consists in studying how to enclose three prescribed volumes in $\R^2$, 
using the least amount of perimeter,
under an additional $1$-dimensional constraint on the intersection of their boundaries, see \eqref{eq:classtr}. 
We will consider $V=\mathcal L^2$ and 
$\mathcal P$ to be the
anisotropic perimeter associated with the so called Grushin plane, defined at the end of the introduction, see \eqref{eq:alphaper}.

We set up the minimal partition problem with trace constraint for a perimeter measure $\mathcal P$ on $\R^2$ and the Lebesgue measure $\mathcal L^2$, in a class of symmetric sets in $\R^2$.
We say that a set $E\subset\R^2$ is {\em $x$-symmetric} (resp. {\em $y$-symmetric})
if $(x,y)\in E$ implies $(-x,y)\in E$
(resp. if $(x,y)\in E$ implies $(x,-y)\in E$).
We say that $E$ is {\em $y$-convex} if the section
$E_x=\{y\in\R : (x,y)\in E\}$ is an interval for every $x\in\R$; finally we say that $E$ is {\em $y$-Schwarz symmetric}
if it is $y$-symmetric and $y$-convex.
We denote by $\S_x$ the class of $\mathcal L^2$-measurable, 
$x$-symmetric sets in $\R^2$ and by $\S^*_y$ the class of
$\mathcal L^2$-measurable and $y$-Schwarz symmetric sets 
in $\R^2$. 

\begin{figure}[h!]
\begin{center}
\vspace{-0.4cm}
\includegraphics[scale=0.5]{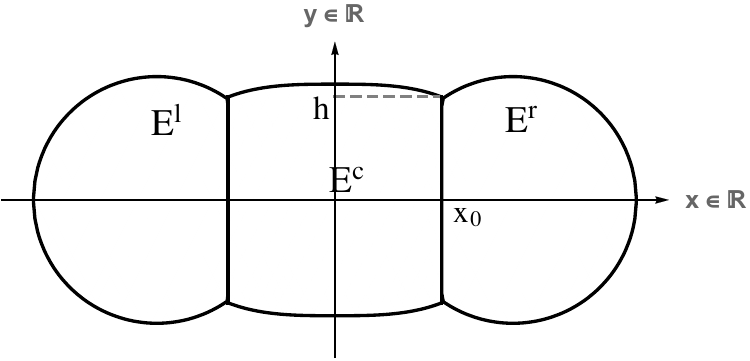}
\caption{A set in the class $\mathcal A$, for $h_1=h_2=h$.}
\label{fig:pmini}
\end{center}
\vspace{-0.4cm}
\end{figure}

Given $v_1,v_2,h_1,h_2\geq 0$, we define the class 
$\mathcal A=\A(v_1,v_2,h_1,h_2)$ of all sets $E\in\S_x\cap\S_y^*$ such that 
for some $x_0>0$, called {\em partitioning point} of $E$, the sets 
\[
E^l=\{(x,y)\in E : x<-x_0\},\quad  E^c=\{(x,y)\in E : |x|<x_0\},\quad E^r=\{(x,y)\in E : x>x_0\}
\] 
satisfy
\vspace{-0.2cm}
\begin{subequations}
\begin{alignat}{2}
\label{eq:classv}
\mathcal L^2(E^c)=v_1,\quad \mathcal L^2(E^l)=\mathcal L^2(E^r)=v_2/2,\\
\label{eq:classtr}[-h_1,h_1]\subset\tr_{x_0-}^xE,\ [-h_2,h_2]\subset\tr_{x_0+}^xE,
\end{alignat}
\end{subequations}
where $\tr_{x_0\pm}^xE$ denote the {\em left} and {\em right traces} of the set $E$ at the point $x_0$, introduced in Definition \ref{def:traces}.
Choosing $h_1=h_2=h>0$, 
the {\em trace constraint} \eqref{eq:classtr} is a relaxed version of the 
equality
\begin{equation}
\label{eq:classtrr}
E_{x_0}
=E_{-x_0}
=[-h,h].
\end{equation}
In other words, a set $E\in\mathcal A$ is such that $E^l,E^c,E^r$ 
have prescribed volumes and $E^c$ 
touches $E^r$ and $E^l$ in segments of a prescribed length, see Figure \ref{fig:pmini}.

We study existence of regular solutions to 
\begin{equation}
\label{eq:pmini}
\inf\{\mathscr P_{\mathcal P}(E) : E\in\mathcal A(v_1,v_2,h_1,h_2)\},\quad  v_1,v_2,h_1,h_2\geq0.
\end{equation}

Our  interest in Problem \eqref{eq:pmini}
comes from the study of the stability of
isoperimetric inequalities.
In the seminal paper \cite{FMP}, the authors 
present a symmetrization technique in the Euclidean space $\R^n$, to prove 
existence of a dimensional constant $C(n)>0$ such that
any $\mathcal L^n$-measurable set $E\subset\R^n$ satisfies 
\begin{equation}
\label{eq:quanti}
P(E)-P(B(0,r_E))\geq C(n)\Big(\min_{x\in\R^n} \mathcal L^n(E\triangle B(x,r_E))\Big)^2.
\end{equation}
Here, $B(0,r)=\{p=(p_1,\dots,p_n)\in\R^n : p_1^2+\dots+p_n^2<r^2\}$, 
and the quantity $r_E\geq0$ is chosen to have $\mathcal L^n(E)=\mathcal L^n(B(0,r_E))$.
Such inequality is known as the sharp {\em quantitative isoperimetric inequality} 
in $\R^n$, see also \cite{CL, FiMP}. 
A minimal partition problem for the Euclidean perimeter $\mathcal P=P$ 
under additional constraints is used in  
\cite{FMP} 
to prove \eqref{eq:quanti} in a class of symmetric sets.
\begin{figure}[h!]
\begin{center}
\vspace{-0.5cm}
\includegraphics[scale=0.4]{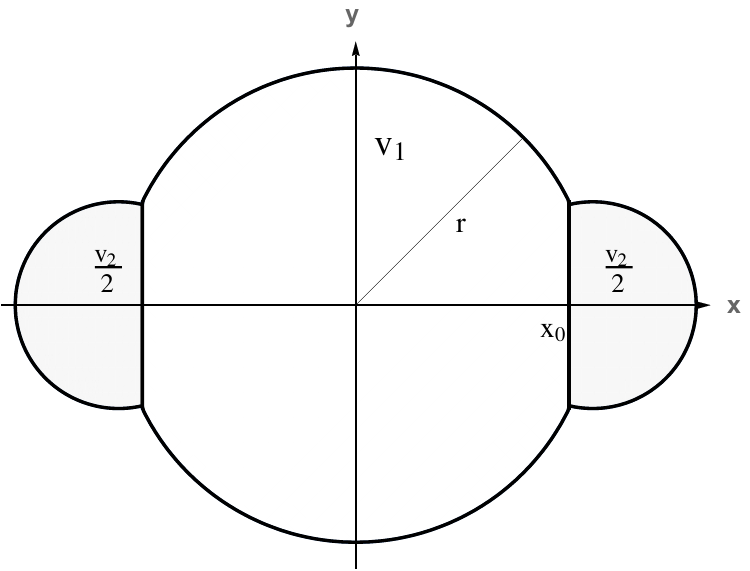}
\vspace{-0.3cm}
\caption{The solution of Problem \eqref{eq:pmini} in the Euclidean setting.}
\label{fig:pmineu}
\end{center}
\vspace{-0.5cm}
\end{figure}

In \cite[Lemma 4.3]{FMP} the authors implicitly use the solution $E$ 
to problem \eqref{eq:pmini} for $\mathcal P=P$, 
with trace constraint given by \eqref{eq:classtrr}.
The solution is, for some $x_0$, $r_0>0$
\[
E=E^l\cup E^c\cup E^r,\quad E^c=B(0,r_0)\cap\{(x,y)\in\R^2 : |x|<x_0\}.
\]
Notice that the central part $E^c$ is the portion of an isoperimetric set lying in a vertical stripe, see Figure \ref{fig:pmineu}.
This is due to the invariance under translations of the standard perimeter  $P$. 

In this paper we study the minimal partition problem with trace constraint in the
{\em Grushin plane}, that is an example of Carnot-Carath\'eodory space, introduced in the context of hypoelliptic operators by Franchi and Lanconelli in \cite{FL82}. 
{\em Carnot-Carath\'eodory spaces} are metric spaces in which the distance is defined in association with a family of vector fields $\mathcal X=\{X_1,\dots,X_r\}$ on a $n$-dimensional manifold, $n\geq r$. The theory of perimeters in such spaces has been developed starting from the 1990s in \cite{CDG, FSSC, GN}, and {\em isoperimetric inequalities} in Carnot-Carath\'eodory spaces are a current object of investigation, see \cite{PRS, P82, GN, FGW, FGuW}. 
Given $\alpha\geq 0$,
the Grushin plane is defined endowing $\R^2$ with the family of vector fields
$\mathcal X_\alpha=\{\partial_x,|x|^\alpha\partial_y\}$,
where $(x,y)$ denotes a point in $\R^2$ and $\partial_x$, $\partial_y$ respectively denote the partial derivative with respect to the first and to the second coordinate.
Given 
a Lebesgue measurable set $E\subset\R^2$, the {\em $\alpha$-perimeter} of $E$
is defined as
\begin{equation}
\label{eq:alphaper}
\P(E)=\sup\Big\{\int_E(\partial_x\varphi_1+|x|^\alpha\partial_y\varphi_2)\;dxdy : \varphi_1,\varphi_2\in C^1_c(\R^2),\ \max_{\R^2}\big(\varphi_1^2+\varphi_2^2\big)^{\frac{1}{2}}\leq1\Big\}.
\end{equation}
If $\mathcal P=P_\alpha$, we set $\mathscr{P}_{\alpha}=\mathscr{P}_{P_\alpha}$, see definition \eqref{eq:mP}. Notice that when $\alpha=0$ the $\alpha$-perimeter is the standard 
De Giorgi perimeter. 
If $E\subset\R^2$ is a bounded set with Lipschitz boundary, 
we have
\begin{equation}
\label{eq:REPR}
\P(E)=\int_{\partial E}\sqrt{(N_x^E(x,y))^2+|x|^{2\alpha}(N_y^E(x,y))^2}\;d\mathcal H^1(x,y),
\end{equation}
where $N^E=(N^E_x,N^E_y)$ is the outer unit normal to  $E$ and $\mathcal H^1$ is the one dimensional Hausdorff measure, see 
\cite[Proposition 2.1]{FM}.
By the representation formula \eqref{eq:REPR}, it is clear that 
$x$-translations
modify the $\alpha$-perimeter with $\alpha>0$, i.e., there exists $E\subset\R^2$ such that $\P(E)=\P(E+(\bar x,0))$, if and only if $|\bar x|=0$.
This is an essential difference with the classical perimeter, or with the one defined in {\em Carnot groups}.

An important feature of the
Grushin plane is that isoperimetric sets are completely characterized.
Given $v>0$, the {\em isoperimetric problem} for the $\alpha$-perimeter is
\begin{equation}
\label{eq:ispb}
\min\{\P(E) : E\subset\R^2,\ \mathcal L^2(E)=v\},
\end{equation}
where $\mathcal L^2$ denotes the two-dimensional Lebesgue measure.
Solutions to \eqref{eq:ispb} have been studied in \cite{MM}, and in \cite{FM} in Grushin structures of dimension $n\geq2$:
up to a vertical translation $\tau_t(x,y)=(x,y+t)$ and an anisotropic dilation $\delta_\lambda(x,y)=(\lambda x,\lambda^{\alpha+1}y)$, the unique solution to problem \eqref{eq:ispb}, called {\em isoperimetric set}, is
\begin{equation}
\label{eq:Ealpha}
\E^\alpha=\{(x,y)\in\R^2 : |y|<\varphi_\alpha(|x|),\ |x|<1\},
\quad \varphi_\alpha(r)=\int_{\arcsin r}^{\frac{\pi}{2}}(\sin{t})^{\alpha+1}\;dt,\ r>0.
\end{equation} 
The one parameter family of dilations $\delta_\lambda$ is such that $\P(\delta_\lambda E)=\lambda^{Q-1}\P(E)$, $\mathcal L^2(\delta_\lambda E)=\lambda^Q\mathcal L^2(E)$, where $Q=\alpha+2$ is called {\em homogeneous dimension}.
In particular, \eqref{eq:Ealpha} implies the validity of the following {\em sharp isoperimetric inequality} for any measurable set $E\subset\R^2$ with finite measure:
\[
\mathcal L^2(E)\leq c(\alpha)\P(E)^{\frac{Q}{Q-1}},\quad c(\alpha)=\frac{\alpha+1}{\alpha+2}\Big(2\int_{0}^\pi\sin^\alpha(t)\;dt\Big)^{-\frac{1}{\alpha+1}}.
\]
When $\alpha=1$, the profile function $\varphi_\alpha$ corresponds to the conjectured isoperimetric profile function of the Heisenberg groups, see \cite{MIsopAx, MR, P82, R, RR1, RR2}.

Quantitative isoperimetric inequalities have been 
studied in Riemannian manifolds providing results in the Gauss space (see \cite{CiaFusMagPra}),
in the $n$-Sphere (see \cite{BDF}) and in the Hyperbolic space (see \cite{BDS}).
Quantitative isoperimetric inequalities in a sub-Riemannian setting are presented 
in \cite{FLM} in the case of the Heisenberg group.
An interesting task would be to study quantitative isoperimetric inequalities in spaces with less isometries, such as the Grushin plane. However, in this paper we will show some  unexpected
obstacles that prevent an adaptation of the techniques in \cite{FMP} to the Grushin plane.
In particular we will show that, when $\alpha>0$, 
solutions to the minimal partition problem \eqref{eq:pmini} are not obtained in their central part as portions of isoperimetric sets lying in a vertical stripe.

From now on we consider $\mathcal P=P_\alpha$ and $\mathscr{P}_\mathcal{P}=\mathscr P_\alpha$. In the first part of this paper we establish existence of solutions.

\begin{thm}
\label{thm:existi}
Let $\alpha\ge0$, $v_1,v_2,h_1,h_2\geq0$. There exists a solution $E=E^l\cup E^c\cup E^r\in\mathcal A(v_1,v_2,h_1,h_2)$ to the minimal partition problem with trace constraint \eqref{eq:pmini} such that $E^c$ is a convex set and $E^l,\ E^r$ have Lipschitz boundaries.
\end{thm}

Minimizers as in Theorem \ref{thm:existi} are called {\em regular}. 
In Proposition \ref{TRpart} we show, under a technical assumption, 
that any regular minimizer $E\in\mathcal A$ 
assumes the least possible traces at the partitioning point $x_0>0$, i.e., $\tr_{x_0-}^x E=[-h_1,h_1]$, and $\tr_{x_0+}^x E=[-h_2,h_2]$.

When $\alpha\in\{0,1\}$ and $v_2=h_2=0$, the geometry of regular solutions can be described in a more precise way.

\begin{thm}
\label{thm:conci}
Let $\alpha\in\{0,1\}$. Given $v_1,h_1> 0$, let
$E=E^l\cup E^c\cup E^r\in\mathcal A(v_1,0,h_1,0)$ be a regular solution for Problem \eqref{eq:pmini} such that 
$\tr_{x_0}^-E=[-h_1,h_1]$. Then 
$E^c=\{(x,y)\in\R^2 : |y|<f(x),\ |x|<x_0\}$, where the function $f$ is given by
\begin{equation}
\label{eq:profi}
f(r)=\lambda^{\alpha+1}\varphi_\alpha\big(\frac{r}{\lambda}\big)+y,
\end{equation}
for some $\lambda=\lambda(\alpha,v_1,h_1)
>0$ and $y
=y(\alpha,v_1,h_1)
\leq 0$ such that  $y=0$ if and only if $\alpha=0$. 
\end{thm}
Due to the presence of the vertical translation $y$ in \eqref{eq:profi}, we deduce by Theorem \ref{thm:conci} that a regular solution of the minimal partition problem
is not obtained as the portion of an isoperimetric set $\delta_\lambda(\E^\alpha)$ lying in 
a vertical stripe, unless $\alpha=0$.
This result shows a delicate point where the techniques of \cite{FMP} fail in the case of the Grushin geometry.

\vspace{0.3cm}

The paper is organized as follows. 
In Section \ref{S:exist}, we prove existence of regular solutions of the minimal partition problem.
The argument is divided into several steps. 
Lemma \ref{lem:1} is an {\em approximation theorem}, that generalizes the classical results in \cite{FSSC}. 
In Lemma \ref{lem:2} we show how to modify a set in the class $\mathcal A$ in order to decrease the perimeter $\mP$ and gain some {\em regularity properties}.
Finally, in Theorem \ref{EXISTpart} we combine {\em lower semicontinuity} of the $\alpha$-perimeter together with a {\em compactness theorem} for sets of finite $\alpha$-perimeter to prove existence of minimizers.

In Section \ref{S:char}, we characterize regular solutions of the minimal partition problem. We call the {\em profile function} of a set $E\in\S_x\cap\S^*_y$, 
the measurable function $f:[0,\infty)\to[0,\infty)$ such that $E=\{(x,y)\in\R^2 : |y|<f(|x|)\}$.
In Propostion \ref{prop:eqdiff} we find differential equations for the profile functionof a regular minimizer.
In Section \ref{SS:traces} we prove, under a technical assumption, that regular solutions of the minimal partition problem \eqref{eq:pmini} satisfy the trace equality \eqref{eq:classtrr}.
We conclude showing formula \eqref{eq:profi} in Proposition \ref{ESTh}, which is proved under the assumptions $\alpha\in\{0,1\}$, $v_1,h_1>0$ and $v_2=h_2=0$. To this purpose, 
we use the differential equations of Proposition \ref{prop:eqdiff} to write the parameters $\lambda$ and $y$ in terms of $\alpha$ and of the given constraints. 

Appendix \ref{app:traces} is dedicated to the notion of trace of a $y$-Schwarz symmetric set.

\section{Existence of minimizers}
\label{S:exist}
In this section we prove existence of solutions to problem \eqref{eq:pmini}.
The proof is divided into several steps that we present in Lemmas \ref{lem:1} and  \ref{lem:2}.

We first introduce some notation. 
For any set $E\subset\R^2$ and $t>0$ we let
\[
\begin{split}
&E_{t-}^x=\{(x,y)\in E : |x|<t\}\quad\text{and}\quad E_{t}^x=\{(x,y)\in E : |x|=t\}\\
\end{split}
\]
In the following, we use the short notation $\{|x|<t\}$ for $\{(x,y)\in\R^2 : |x|<t\}$.

%

The next lemma is a refinement of the approximation theorem by smooth sets for sets with finite $\alpha$-perimeter, see \cite[Theorem 2.2.2]{FSSC}. We say that  a set $E\subset\R^2$ is {\em locally Lipschitz (resp. locally $C^\infty$)} if its boundary $\partial E$ is a locally Lipschitz (resp. locally $C^\infty$) curve.

\begin{lem}[Approximation by smooth sets]
\label{lem:1}
Given $v_1,v_2,h_1,h_2>0$, let $E\in\mathcal A(v_1,v_2,h_1,h_2)$ be a set with finite $\alpha$-perimeter. 
Let $x_0>0$ be the partitioning point for $E$ and $y_0^\pm$ be such that $\tr_{x_0\pm}^xE=[-y_0^\pm,y_0^\pm]$. 
Then there exists a sequence of locally $C^\infty$ sets $\mathcal E_j \in\S_x\cap\S^*_y$, $j\in\N$
such that
\begin{subequations}
\begin{alignat}{3}
\label{eq:lem11}\bullet&\lim_{j\to\infty}\P( (\mathcal E_j)_{x_0-}^x)=\P(E_{x_0-}^x)\text{ and }
\lim_{j\to\infty}\P( \mathcal E_j\setminus (\mathcal E_j)_{x_0-}^x)
=\P(E\setminus E_{x_0-}^x);\\
\label{eq:lem12}\bullet&\lim_{j\to\infty}\mathcal L^2( (\mathcal E_j)_{x_0-}^x)=\mathcal L^2(E_{x_0-}^x)\text{ and }
\lim_{j\to\infty}\mathcal L^2(\mathcal E_j\setminus (\mathcal E_j)_{x_0-}^x)=\mathcal L^2( E\setminus E_{x_0-}^x);\\
\label{eq:lem13}\bullet&\text{ if $\tr_{x_0\pm}^x \mathcal E_j=[-q_j^\pm,q_j^\pm]$, for some $q_j^\pm\geq0$, we have }
q_j^\pm\to y_0^\pm\text{ as }j\to\infty.
\end{alignat}
\end{subequations}
\end{lem}

\proof

Let 
$J\in C^\infty_c(B(0,1))$, be a positive symmetric mollifier and, for 
$\varepsilon>0$,
let $J_\varepsilon(p)=\frac{1}{\varepsilon^2}J(|p|/\varepsilon)$, 
$p\in\R^2$.
Following \cite[Theorem 1.24]{G},
we can choose $t\in(0,1)$
such that, given any sequence $\varepsilon_j\to0^+$, $j\to\infty$
the sets
$\mathcal E_j=\{p\in\R^2 : J_{\varepsilon_j}*\chi_{E}(p)>t\}$
are locally smooth sets satisfying \eqref{eq:lem12}, and
\begin{equation}
\label{vrperim}
\begin{split}
&\lim_{j\to\infty}\P(\mathcal E_j;(\mathcal E_j)_{x_0-}^x)
=\P(E; E_{x_0-}^x),\
\lim_{j\to\infty}\P(\mathcal E_j;\mathcal E_j\setminus (\mathcal E_j)_{x_0-}^x\})=\P(E;E\setminus E_{x_0-}^x).
\end{split}
\end{equation}

By symmetry of $J$, the sets $\mathcal E_j$ are $x$-symmetric.
We show that they are also $y$-Schwarz
symmetric. In fact, since $E$ is $y$-Schwarz symmetric,
we have $\chi_E(\bar x-x',\bar y-y')\leq\chi_E(\bar x-x',y-y')$ for $(x',y')\in\R^2$ and $|y|<|\bar y|$. 
Hence, setting $h_{\varepsilon_j}=J_\varepsilon*\chi_E$, we obtain
\[\begin{split}
t<h_{\varepsilon_j}(\bar x,\bar y)=\int_{B_{\varepsilon_j}(0)}&\!\!\!\!\!\!J_\varepsilon(x',y')\chi_E(\bar x-x',\bar y-y')\;dx'dy'\\
&\leq\int_{B_{\varepsilon_j}(0)}\!\!\!\!\!\!J_\varepsilon(x',y')\chi_E(\bar x-x',y-y')\;dx'dy'=h_{\varepsilon_j}(\bar x,y),\end{split}
\]
which implies $(\bar x,y)\in\mathcal E_j$ for $|y|<|\bar y|$.

We prove claim \eqref{eq:lem13}.
Let $\phi_j$ denote the profile function of $\mathcal E_j$.
Then we have $\tr_{x_0\pm}^x\mathcal E_j=[-q_j^\pm,q_j^\pm]$, where (see Remark \ref{tracesm})
\[
q_j^-=\lim_{x\to x_0^-}\phi_j(x_0),\quad q_j^+=\lim_{x\to x_0^+}\phi_j(x_0).
\]
In the same way, for $0<\sigma< y_0^+$, there exists $\delta=\delta(\sigma)>0$
such that 
\begin{equation}
\label{ppicc}
|f(x)- y_0^+|<\sigma \quad\text{for}\quad x_0< x< x_0+\delta,
\end{equation}
where $f$ is the profile function of $E$.
We first claim that  there exists $\bar j=\bar j(\sigma)$ such that
\begin{equation}
\label{ppicci}
(x-\xi, y-\eta)\in E\quad\text{for}\quad (\xi,\eta)\in B(0,\varepsilon_j)
\end{equation} 
for any $j\geq\bar{j}$ and for $(x,y)$ sufficiently far from $(x_0,y_0^+)$.
To prove it, choose $\bar j(\sigma)\in\N$ to have $\varepsilon_{\bar j}<\min\{\sigma,\ \delta(\sigma)/4\}$ and let $(x,y)\in A_{\sigma}=(x_0+\varepsilon_j,x_0+\frac{\delta}{2})\times(0,y_0^+-2\sigma)$ for $j\geq \bar j$. Then, 
for $-\varepsilon_j<\xi,\eta<\varepsilon_j$,
we get $x-\xi\in(x_0,x_0+\delta)$. Hence \eqref{ppicc} leads to 
$y-\eta<f(x-\xi)$ and the claim is proved.

We deduce from \eqref{ppicci} that $A_\sigma\subset \mathcal E_j$ for $j\geq j(\sigma).
$
In fact,
if $(x,y)\in A_\sigma$ we have
\[
h_{\varepsilon_j}(x,y)
=\int_{B(0,\varepsilon_j)}J_{\varepsilon_j}(\xi,\eta)\chi_E(x-\xi,y-\eta)\;d\xi d\eta
=\int_{B(0,\varepsilon_j)}J_{\varepsilon_j}(\xi,\eta)\;d\xi d\eta=1>t.
\]
In particular,
\[
(-y_0^++2\sigma,y_0^+-2\sigma)\subset \tr_{(x_0+\varepsilon_j)+}^x\mathcal E_j\quad\text{for every}\quad j>\bar j(\sigma).
\]
Similarly,
we can choose $\bar{\bar j}(\sigma)\in\N$
such that 
\[
\tr_{(x_0+\varepsilon_j)+}\mathcal E_j
\subset(-y_0^+-2\sigma,y_0+2\sigma)\quad\text{for}\quad j\geq\bar{\bar{j}}(\sigma),
\]
and \eqref{eq:lem13} for the right trace is proved. The argument to prove \eqref{eq:lem13} for the left trace is analogous. Statement \eqref{eq:lem11} follows from \eqref{eq:lem13} and \eqref{vrperim}.
\endproof


\begin{lem}[Regularization]
\label{lem:2}
Let $v_1,v_2,h_1,h_2\geq 0$ and $E\in\mathcal A(v_1,v_2,h_1,h_2)$ be a locally $C^\infty$-set with finite $\alpha$-perimeter. Then, 
there exists a set $\tilde E\in\mathcal A(v_1,v_2,h_1,h_2)$ such that, if $\tilde x_0$ is the partitioning point for $\tilde E$, there holds:
\begin{enumerate}
\item $\tilde E_{\tilde x_0-}^x$ is convex and $\tilde E\setminus\tilde E_{\tilde x_0-}^x$ has locally Lipschitz boundary;
\item $\mathscr P_\alpha(\tilde E)\leq\mathscr P_\alpha(E)$, in particular $\P(\tilde E_{\tilde x_0-}^x)\leq \P(E_{x_0-}^x)$ and $\P(\tilde E\setminus \tilde E_{\tilde x_0-}^x)\leq \P(E\setminus E_{x_0-}^x)$ where $x_0$ is the partitioning point for $E$.
\end{enumerate}
\end{lem}

\proof 
Let 
$\mathrm{tr}_{x_0-}^xE=[-q^-,q^-]$ and $\mathrm{tr}_{x_0+}^xE=[-q^+,q^+]$.
We divide the proof into the following steps, corresponding to operations performed on the set $E$.
\vspace{0.2cm}

\noindent\textit{Step 1. (Gluing around the $y$-axis)}.
Starting from $E$, we construct a set
$\hat E\in\S_x\cap\S^*_y$
such that
there exist $0<\hat x_0\leq x_0$
satisfying:
\begin{enumerate}
\item the Euclidean outer unit normal to $\hat E$ exists outside a set of $\mathcal H^1$-measure zero;
\item if $\hat\phi:[0,\infty)\to[0,\infty)$ denotes the profile function of $\hat E$ and $\hat D=\inf\{d\geq0 : \hat\phi(x)=0\text{ for }x\geq d\}$, then
$\displaystyle \mathcal H^1(\{x\in[0,\hat D] : \hat\phi(x)=0\})=0$;
\item $\displaystyle \P(\hat E_{\hat x_0-}^x)\leq\P(E_{x_0-}^x)$
and $\displaystyle \P(\hat E\setminus \hat E_{\hat x_0-}^x)\leq\P(E\setminus E_{x_0-}^x)$;
\item $\displaystyle \mathcal L^2(\hat E_{\hat x_0-}^x)=\mathcal L^2(E_{x_0-}^x)$ and
 $\displaystyle \mathcal L^2(\hat E\setminus\hat E_{\hat x_0-}^x)=\mathcal L^2( E\setminus E_{x_0-}^x)$;
\item $\displaystyle \tr_{\hat x_0-}^x\hat E=\tr_{x_0}^{-}  E$
and $\displaystyle \tr_{\hat x_0+}^x\hat E=\tr_{x_0}^{+}  E$.
\end{enumerate}
\begin{figure}[h!!!]
\centering
\includegraphics[scale=0.4]{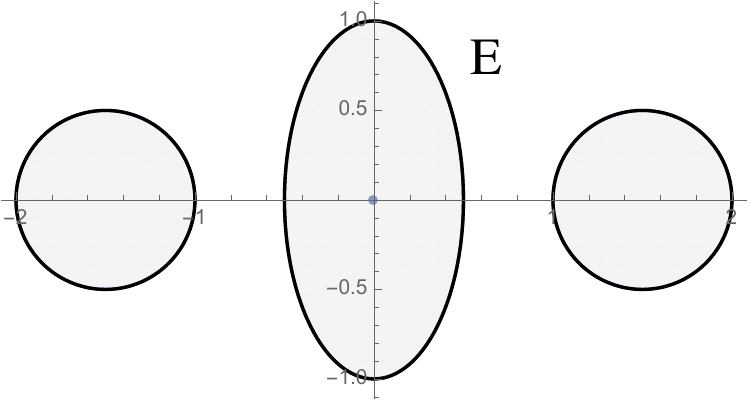} 
\includegraphics[scale=0.4]{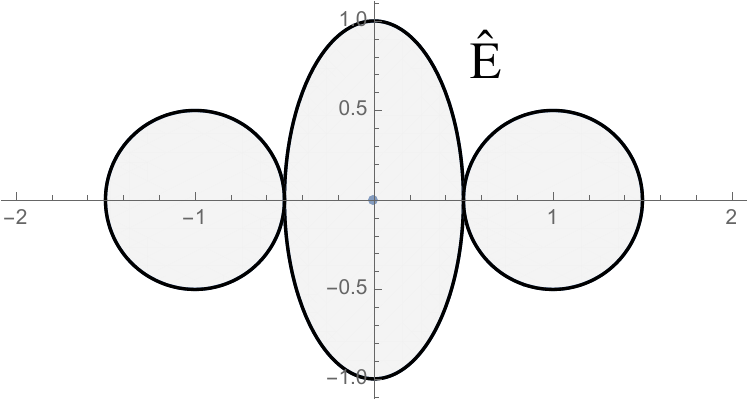}
\caption{From $ E$ to $\hat E$.}
\label{p_holes}
\end{figure}

Let $\phi:[0,\infty)\to[0,\infty)$ be the profile function of $E$.
Let $D=\inf\{d\geq0 : \phi(x)=0\text{ for }x\geq d\}$. 
Define the set
$Z:=\{x\in[0,D] :  \phi(x)=0\}$ and write
$Z=Z^1\cup Z^2$ with
\[
\begin{split}
Z^1&=\{x\in[0,D] : \phi(x)=0 \text{ and }  \phi(\xi)\neq0\text{ for }
\xi\in(x-\delta,x+\delta)\setminus\{x\}\text{ for some }\delta>0\},\\
Z^2&=\{x\in[0,D] : \exists\delta>0 : \phi(\xi)=0
\text{ for }\xi\in(x-\delta,x]\text{ or }\xi\in[x,x+\delta)\}.
\end{split}
\]
By symmetry and smoothness of $ E$, 
we have $Z^1=\emptyset$.
On the other hand, the set $Z^2$ is the union of at most 
countably many open intervals, 
being the complement
in $\R$ of the closed set $\{x\in\R : (x,0)\in\overline{ E}\}$.
Let
\[\begin{split}
&Z^2=\bigcup_{k\in \mathfrak I}(a_k,b_k)
\cup \bigcup_{k\in \mathfrak J}(c_k,d_k)
\cup \bigcup_{k\in \mathfrak I}(-b_k,-a_k)
\cup \bigcup_{k\in \mathfrak J}(-d_k,-c_k),
\end{split}\]
for $0\leq a_1<b_1<a_2<b_2<\dots\leq x_0<c_1<d_1<c_2<d_2<\dots\leq D$
and $\mathfrak I,\mathfrak J\subset\N$.
We rearrange $ E$ 
in at most countably many steps,
each one corresponding to
an interval $(a_k,b_k)$. 

Associated with the interval $(a_1,b_1)$, define the set
\[\begin{split}
 E_1&
= ( E)_{a_1-}^x 
\cup\{(x+a_1-b_1,y) : (x,y)\in  E,\ x> b_1\}
\cup\{(x+b_1-a_1,y) : (x,y)\in  E,\ x<- b_1\}
\end{split}\]
obtained by ``filling the hole" $(a_1,b_1)$. $E$ is $x$-symmetric and $y$-Schwarz symmetric. Moreover, letting $x_1=x_0+a_1-b_1<x_0$, we have
\[
\mathcal L^2(( E_1)_{x_1-}^x)
=\mathcal L^2( E_{x_0-}^x),\quad
\mathcal L^2( E_1\setminus( E_1)_{x_1-}^x)
= \mathcal L^2( E\setminus E_{x_0-}^x)
\]
and
$\tr_{x_1\pm}^x E_1=\tr_{x_0\pm}^x E$.
We show that $\P\big((E_1)_{x_1-}^x\big)\leq\P\big(E_{x_0-}^x\big)$. 
Let $N_1(p)=(N_{1x}(p),N_{1y}(p))$
be the Euclidean outer unit normal to $ E_1$ 
at $p=(x,y)\in\R^2$, for $|x|\neq a_1$.
If $N=(N_x,N_y)$ is the Euclidean outer unit normal
to $\partial  E$,
we have
\[\begin{split}
&N_1(x,y)=N(x-a_1+b_1,y)
\text{ if }  x> a_1
\ \text{ and }\ 
N_1(x,y)=N(x-b_1+a_1,y)
\text{ if } x< -a_1.
\end{split}
\]
The claim then follows by the representation formula \eqref{eq:REPR},
observing that, given $x\in \partial\big(( E_1)_{x_1-}^x\setminus (E_1)_{a_1-}^x\big)$
we have
\[
|x|^{2\alpha}N_{y}(x-a_1+b_1,y)^2= |\bar x-b_1+a_1|^{2\alpha}{N_{y}(\bar x, y)}^2
\leq |\bar x|^{2\alpha}N_y(\bar x,y)^2,
\]
for $\bar x=x-a_1+b_1\in\partial\big(E_{x_0}^x\setminus E_{b_1}^x\big)$.


We repeat this procedure for every interval $(a_k,b_k)$, $k\in\mathfrak I$. 
After at most countably many operations,
we obtain a $x$-symmetric and $y$-Schwarz symmetric
set $\hat E$ satisfying claims 3, 4 and 5 for
$
\hat x_0=x_0-\sum_{i\in\mathcal I}(b_i-a_i).
$
Repeating this argument for the intervals $(c_k,d_k)$, $k\in\mathfrak J$,
we define the set $\hat E$ that satisfies also claims
1 and 2. In particular, letting
\[
\hat{Z}=\Big\{a_k-\sum_{i=1}^{k-1}(b_i-a_i) : k\in\mathfrak I\Big\}\cup
\{c_k-\sum_{i=1}^{k-1}(d_i-c_i) : k\in\mathfrak J\},
\]
the outer unit normal to $\hat E$ exists outside
the set $\{(x,y)\in\R^2 : \ |x|\in \hat Z\}$.

\vspace{0.2cm}

\noindent\textit{Step 2. (Reflection in the vertical direction)}
We rearrange the set $\hat E$ into a
$x$-symmetric and $y$-Schwarz symmetric set
$\hat{\hat{E}}$ with profile function 
$\hat{\hat{\phi}}:[0,\infty)\to[0,\infty)$
such that
\begin{enumerate}
\item The Euclidean outer unit normal to $\hat{\hat{E}}$ exists outside
a set of $\mathcal H^1$-measure zero;
\item $\displaystyle \hat{\hat{\phi}}(|x|)\geq q^-$ for $x\in\R$, $|x|<\hat x_0$;
\item $\displaystyle \P((\hat{\hat{E}})_{\hat x_0-}^x)\leq\P(\hat E_{\hat x_0-}^x)$
and $\displaystyle \P(\hat{\hat{E}}\setminus (\hat{\hat{E}})_{\hat x_0-}^x)=\P(\hat E\setminus(\hat E)_{\hat x_0-}^x)$;
\item $\displaystyle 
\mathcal L^2((\hat{\hat{E}})_{\hat x_0-}^x)\geq 
\mathcal L^2(\hat E_{\hat x_0-}^x)$
and
$\displaystyle 
\mathcal L^2(\hat{\hat{E}}\setminus(\hat{\hat{E}})_{\hat x_0-}^x)=
\mathcal L^2(\hat{E}\setminus(\hat{E})_{\hat x_0-}^x)$;
\item $\displaystyle \tr_{\hat x_0-}^x\hat{\hat{E}}
=\tr_{\hat x_0-}^x\hat{E}$ and $\tr_{\hat x_0+}^x\hat{\hat{E}}=\tr_{\hat x_0+}^x\hat{E}$.
\end{enumerate}

\begin{figure}[h!!!]
\centering
\includegraphics[scale=0.35]{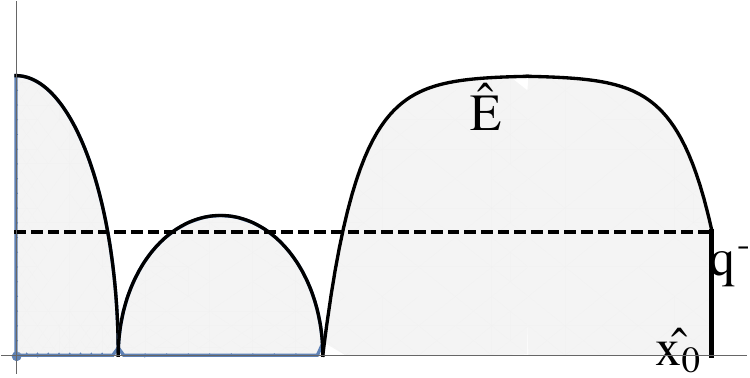} 
\hspace{0.3cm}
\includegraphics[scale=0.35]{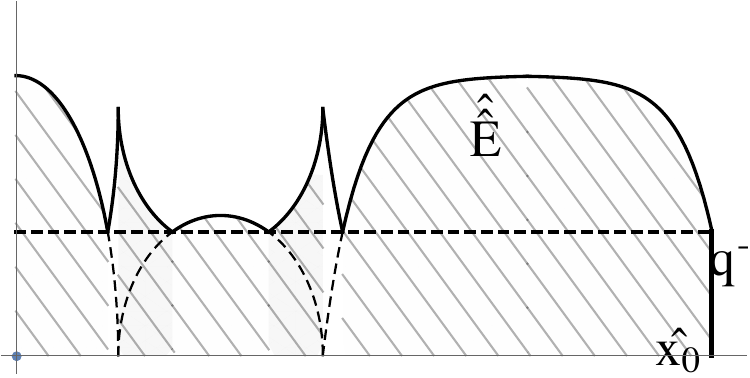}
\caption{The set $\hat E$ and the rearranged $\hat{\hat{E}}$.}
\label{p_ACCA}
\end{figure}

We define the rearranged function $\hat{\hat \phi}:[0,\infty)\to[0,\infty)$,
\[\hat{\hat{\phi}}(x)=
\left\{
\begin{array}{lll}
&|\hat{\phi}(x)-q^-|+q^-
& \text{ if }|x|<\hat x_0, \\
&\hat{\phi}(x)&\text{ if }|x|>\hat x_0.
\end{array}
\right.
\]
Let $\hat{\hat{E}}$ be the $x$- and $y$-symmetric set generated
by $\hat{\hat{\phi}}$ (see Figure \ref{p_ACCA}). 

Claims 2, 4 and 5 are obvious. The points where this operation breaks the regularity of the boundary are the ones in 
$\mathcal K= 
\{(x,\hat \phi(x))\in\partial \hat{E} : 
\hat \phi(x)=q^-\text{ and }\exists\delta>0 : \hat \phi(\xi)\neq q^-\text{ for }
\xi\in(x,x+\delta)\text{ or }\xi\in(x-\delta,x)\}$. Claim 1 follows since $\mathcal K$ is countable.

Moreover, letting $\hat{\hat{N}} =(\hat{\hat{N}}_{x},\hat{\hat{N}}_{y})$ be the 
outer unit
normal to $\hat{\hat{ E}}$ outside $\mathcal K\cup \hat Z$
and 
$\hat N=(\hat N_{x},\hat N_{y})$ be the outer unit normal to 
$\partial\hat E$, we have for any $(x,y)\in\partial \hat E\setminus 
\big(\mathcal K \cup \hat Z\big)$,
$|x|<\hat x_0$,
\[
\hat{\hat N} (x,|y-q^-|+q^-)=\big(\hat N_{x}(x,y),\mathrm{sgn}(y-q^-)\hat N_{y}(x,y)\big).
\]
We hence deduce claim 3 by the Representation formula \eqref{eq:REPR}, by writing
\[\begin{split}
\P\big((\hat{\hat{E}})_{\hat x_0-}^x\big)
&=\int_{\{p\in\partial\hat E_{\hat x_0-}^x :\ \hat N_{y}(p)\neq 0\}}
\sqrt{\hat N_{x}^2+|x|^{2\alpha}\hat N_{y}^2}\;d\mathcal H^1
+\mathcal H^1\big(\{p\in\partial (\hat{\hat{E}})_{\hat x_0-}^x :\ \hat {\hat{N}}_{y}(p)= 0\}\big).
\end{split}
\]

\textit{Step 3. (Convexification and regularization)}
We finally rearrange the set $\hat{\hat{E}}$
into a $x$-symmetric and $y$-Schwarz symmetric set
$\tilde{E}$, such that there exists
$0<\tilde x_0<x_0$ satisfying:
\begin{enumerate}
\item $\displaystyle (\tilde{E})_{\tilde x_0-}^x$ is convex, 
and $\tilde E\setminus\tilde E_{\tilde x_0-}^x$
is locally lipschitz;
\item $\displaystyle \P((\tilde{E})_{\tilde x_0-}^x)\leq \P( E_{x_0-}^x)$ and 
$\displaystyle \P(\tilde{E}\setminus (\tilde{E})_{\tilde x_0-}^x)\leq \P( E\setminus  E_{x_0-}^x)$;
\item $\displaystyle \mathcal L^2((\tilde{E})_{\tilde x_0-}^x)=\mathcal L^2( E_{x_0-}^x)$ and
$\displaystyle \mathcal L^2(\tilde{E}\setminus(\tilde{E})_{\tilde x_0-}^x)=\mathcal L^2( E\setminus E_{x_0-}^x)$;
\item $\tr_{\tilde x_0-}^x\tilde{E}\supset [-q^-,q^-]$ and 
$\tr_{\tilde x_0+}^x\tilde{E}\supset [-q^+,q^+]$.
\end{enumerate}
This will conclude the proof.

We introduce the functions
\[
\Psi,\Phi:\R^2\to\R^2,\quad \Psi(x,y)=\big(\mathrm{sgn}(x)\frac{|x|^{\alpha+1}}{\alpha+1},y\big),
\quad \Phi(\xi,\eta)=\big(\mathrm{sgn}(\xi)|(\alpha+1)\xi|^{\frac{1}{\alpha+1}},\eta\big),
\]
which are homeomorphisms such that $\Psi^{-1}=\Phi$.
As shown in \cite[Prop. 2.3]{MM}, for any measurable set $F\subset\R^2$, we have
\[
\P(F)=P(\Psi(F)) \quad\text{and}\quad \mathcal L^2(F)=\mu(\Psi(F)),
\] 
where $P$ denotes the Euclidean perimeter and
$\mu$ is a Borel measure on $\R^2$ defined on Borel sets as follows:
\[
\mu(A)=\int_A|(\alpha+1)\xi|^{-\frac{\alpha}{\alpha+1}}\;d\xi d\eta,\quad A\subset\R^2\text{ Borel}.
\]

Let $F^{\mathrm{c}}=\Psi((\hat{\hat{E}})_{\hat x_0-}^x)\subset\R^2$
and consider its convex envelope in $\R^2$, $\mathrm{co}(F^{\mathrm{c}})$. 
Since $\Phi$, $\Psi$ preserve the
symmetries, the transformed
set $\mathcal E^{\mathrm{c}}=\Phi(\mathrm{co}(F^{\mathrm{c}}))$ is $x$- and $y$-Schwarz symmetric. 
Arguing as in \cite[p. 362]{MM} we obtain that $\mathcal E^{\mathrm{c}}$ is also convex.
Moreover, we have
\begin{subequations}
\begin{alignat}{3}
\label{volprt}
&\mathcal L^2(\mathcal E^{\mathrm{c}})=\mu(\mathrm{co}(F^{\mathrm{c}}))\geq \mu( F^{\mathrm{c}})
=\mathcal L^2((\hat{\hat{E}})_{\hat x_0-}^x)
\geq \mathcal L^2(\hat E_{\hat x_0-}^x)
\geq \mathcal L^2( E_{x_0-}^x), \\
\label{perprt}
&\P(\mathcal E^{\mathrm{c}})=P(\mathrm{co}(F^{\mathrm{c}}))\leq P(F^{\mathrm{c}})= \P((\hat{\hat{E}})_{\hat x_0-}^x),\\
\label{trprt}
&\tr_{\hat x_0-}^x\mathcal E^{\mathrm{c}}=[-q^-,q^-].
\end{alignat}
\end{subequations}
By \eqref{volprt},
we define $\tilde x_0\in[0,\hat x_0]$
such that
$\mathcal L^2\big((\mathcal E^{\mathrm{c}})_{\tilde x_0-}^x\big)
=\mathcal L^2\big( E_{x_0-}^x\big)$.
Using a calibration argument as in \cite[Proposition 4.2]{FM}, we deduce 
that $P_\alpha((\mathcal E^{\mathrm{c}})_{\tilde x_0-}^x)\leq P_\alpha(\mathcal E^{\mathrm{c}})$.
Moreover, $[-q^-,q^-]\subset
\tr_{\tilde x_0-}^x\mathcal E^{\mathrm{c}}$ since $\mathcal E^{\mathrm{c}}$ has a monotone decreasing profile function.
\begin{figure}[h!]
\begin{center}
\includegraphics[scale=0.45]{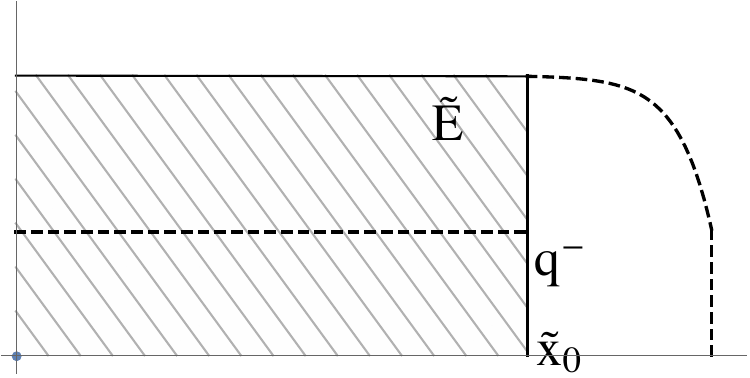}
\caption{The convexified set $F^c$, cut at $\tilde x_0$.}
\end{center}
\end{figure}

Define the set
$ \mathcal E =(\mathcal E^{\mathrm{c}})_{\tilde x_0-}^x
\cup
\{(x-\hat x_0+\tilde x_0,y)\in\R^2 : (x,y)\in\hat{\hat{E}},\ x>\hat x_0\}
\cup
\{(x+\hat x_0-\tilde x_0,y)\in\R^2 : (x,y)\in\hat{\hat{E}},\ x<-\hat x_0\}$.
Arguing as in Step 1, we have $\P( \mathcal E \setminus( \mathcal E )_{\tilde x_0-}^x)
\leq \P(\hat{\hat{E}}\setminus(\hat{\hat{E}})_{\hat x_0-}^x)$,
$\mathcal L^2( \mathcal E \setminus( \mathcal E )_{\tilde x_0-}^x)
=\mathcal L^2(\hat{\hat{E}}\setminus(\hat{\hat{E}})_{\hat x_0-}^x)$
and $\tr_{\tilde x_0+}^x \mathcal E =\tr_{\hat x_0+}^x\hat{\hat{E}}$.

Now, the same argument used to prove \eqref{volprt}-\eqref{trprt}, shows that the sets
$ \mathcal E ^r=\Phi(\mathrm{co}(\Psi( \mathcal E \cap\{x>\tilde x_0\})))$ and 
$ \mathcal E ^l=\{(-x,y) : (x,y)\in  \mathcal E ^r\}$
satisfy
\[\begin{split}
&\mathcal L^2( \mathcal E ^r\cup  \mathcal E ^l)\geq \mathcal L^2( \mathcal E \setminus ( \mathcal E )_{\tilde x_0-}^x),
\quad\tr_{\tilde x_0+} \mathcal E ^r=\tr_{\tilde x_0+}  \mathcal E ,\quad
\P( \mathcal E ^l\cup \mathcal E ^r)\leq \P( \mathcal E \setminus ( \mathcal E )_{\tilde x_0-}^x).
\end{split}
\]
Moreover, since $\mathrm{co}(\Psi( \mathcal E \cap\{x>\tilde x_0\}))$ is a convex set and $\Phi$ is $C^{\infty}$-smooth outside 0,
\textcolor{black}{the set $ \mathcal E ^r$ is locally Lipschitz.}
Let $r\geq \tilde x_0$ be such that 
$\mathcal L^2(( \mathcal E ^r\cup  \mathcal E ^l)_{r-}^x)
=\mathcal L^2( E\setminus  E_{x_0-}^x)$.
The set
\begin{equation}
\label{errej}
\tilde{E}=
( \mathcal E ^c)_{\tilde x_0-}^x\cup
( \mathcal E ^r\cup  \mathcal E ^l)_{r-}^x
\end{equation}
satisfies all the claims of this Step. This concludes the proof.
\endproof

\begin{rem}
\label{rem:lip} 
Given a set $E\in\mathcal A(v_1,v_2,h_1,h_2)$ with partitioning point $x_0>0$ it is always possible to construct a set $\tilde E\in\mathcal A(v_1,v_2,h_1,h_2)$ such that 
$\mathscr P_\alpha(\tilde E)\leq \mathscr P_\alpha(E)$, $\tilde E_{x_0-}^x=E_{x_0-}^x$, and $\tilde E\setminus\tilde E_{x_0-}^x$  is locally Lipschitz. This can be done in the same way as to construct 
the set $\tilde E$ defined in \eqref{errej}.
\end{rem}

We prove that Problem \eqref{eq:pmini} for $\P$ admits regular solutions.
Given $v_1,v_2,h_1,h_2\geq 0$, let
\[
C_{MP}=\inf\{\mathscr P_\alpha(E) : E\in\mathcal A(v_1,v_2,h_1,h_2)\}.
\]
The constant $C_{MP}$ is positive thanks to the validity of the following isoperimetric inequality:
for any $\mathcal L^2$-measurable set $E\subset\R^2$ with finite measure
 \[
\P(E)\geq C\mathcal L^2(E)^{\frac{\alpha+1}{\alpha+2}}
\]
for some geometric constant $C>0$, see
\cite{FGW, FGuW, GN} (see also \cite[Proposition 1.3.4]{FrThesis}).
In fact, by
the formulas
\[
\P(E_{t-}^x)=\P(E;E_{t-}^x)+\mathcal H^1(E_t^x)\quad\text{and}\quad
\P(E\setminus E_{t-}^x)=\P(E;E\setminus E_{t-}^x)+\mathcal H^1(E_t^x),
\]
holding for a.e. $t>0$ and any set $E$ with finite measure and finite $\alpha$-perimeter (see \cite[Proposition 4.1]{FM}), we deduce that for $E\in\mathcal A$, $\mP(E)\geq P(E)$. 
Equality holds if and only if $E\in\mathcal A(v_1,v_2,h,h)$ is
such that $\tr_{x_0}^{-}E=\tr_{x_0}^{+}E=[-h,h]$, for $h>0$.

\begin{thm}\label{EXISTpart}
Let $v_1,v_2,h_1,h_2\geq 0$. There exists a bounded set $E\in\mathcal A$
with partitioning point $x_0\geq 0$ realizing the infimum in \eqref{eq:pmini} for $\mathscr P_\alpha$ and such that
$E_{x_0-}^x$ is convex, {\color{black} and 
$E\setminus E_{x_0-}^x$
is locally Lipschitz.}
\end{thm}

\proof
Let $(E_m)_{m\in\N}$ be a minimizing 
sequence for the infimum in \eqref{eq:pmini}, namely
\[
E_m\in\mathcal A\quad \mathscr P_\alpha(E_m)\leq C_{MP}\Big(1+\frac{1}{m}\Big)\quad m\in\N.
\]
Let $x_m>0$ be the partitioning point for $E_m$ and 
$f_m:[0,\infty)\to[0,\infty)$  be its profile function.
Moreover, let $y_m^-,y_m^+\geq0$ be such that
$\tr_{x_m\pm}^xE_m=[-y_m^\pm, y_m^\pm]$.
By Remark \ref{tracesm},
\[
\lim_{x\to x_m^+}f_m(x)=y_m^+,\quad \lim_{x\to x_m^-}f_m(x)=y_m^-,
\quad\text{with}\quad
y_m^+\geq h_2,\ 
y_m^-\geq h_1.
\]

Let $m\in\N$. By Lemma \ref{lem:1}, let $\mathcal E_j^m\in\S_x\cap\S^*_y$, $j\in\N$ be a sequence of smooth sets satisfying \eqref{eq:lem11}-\eqref{eq:lem13} for $E=E_m$. 
Define $q_{jm}^\pm\geq 0$ such that $\mathrm{tr}_{x_m\pm}^x\mathcal E_j^m=[-q_{jm}^\pm,q_{jm}^\pm]$.
For any $j\in\N$, apply Lemma \ref{lem:2} to $\mathcal E_j^m$, obtaining 
a $x$-symmetric and $y$-Schwarz symmetric set
$\tilde{\mathcal E}_j^m$, such that there exists
$0<\tilde x_j^m<x_m$ satisfying:
\begin{enumerate}
\item $\displaystyle (\tilde{\mathcal E}_j^m)_{\tilde x_j^m-}^x$ is convex, 
and 
$\tilde{\mathcal E}_j^m\setminus (\tilde{\mathcal E_j^m})_{\tilde x_j^m-}^x$ is locally lipschitz;
\item $\displaystyle \P((\tilde{\mathcal E}_j^m)_{\tilde x_j^m-}^x)\leq \P((\mathcal E_j^m)_{x_m-}^x)$ and 
$\displaystyle \P(\tilde {\mathcal E}_j^m\setminus (\tilde{\mathcal E}_j^m)_{\tilde x_j^m-}^x)\leq \P(\mathcal E_j\setminus (\mathcal E_j^m)_{x_m-}^x)$;
\item $\displaystyle \mathcal L^2((\tilde{\mathcal E}_j)_{\tilde x_j-}^x)=\mathcal L^2((\mathcal E_j)_{x_0-}^x)$ and
$\displaystyle \mathcal L^2(\tilde{\mathcal E}_j\setminus(\tilde{\mathcal E}_j)_{\tilde x_j-}^x)=\mathcal L^2(\mathcal E_j\setminus(\mathcal E_j)_{x_0-}^x)$;
\item $\tr_{\tilde x_j^m-}^x\tilde{\mathcal E}_j^m\supset [-q_j^-,q_j^-]$
and $\tr_{\tilde x_j+}^x\tilde{\mathcal E}_j\supset [-q_j^+,q_j^+]$.
\end{enumerate}

By \eqref{eq:lem13}, for any $m\in\N$ there exists 
$J(m)\in\N$ such that for $j\geq J(m)$,
we have 
\begin{equation}
\label{trac}
|q_{jm}^\pm-y_m^\pm|\leq \frac{1}{m}
\end{equation}
and
\begin{equation}
\label{feder}
\begin{split}
&\Big|\P\big((\tilde {\mathcal E}_j^{m})_{\tilde x_j^m-}^x\big)-\P((E_m)_{x_m-}^x)\big)\Big|\leq \frac{1}{m},
\quad
\Big|\P\big(\tilde {\mathcal E}_j^{m}\setminus (\tilde {\mathcal E}_j^{m})_{\tilde x_j^m-}^x\big)-\P\big(E_m\setminus (E_m)_{x_m-}^x\big)\Big|\leq \frac{1}{m},\\
&\Big|\mathcal L^2\big((\tilde {\mathcal E}_j^{m})_{\tilde x_j^m-}^x\big)-\mathcal L^2((E_m)_{x_m-}^x)\big)\Big|\leq \frac{1}{m},
\quad
\Big|\mathcal L^2\big(\tilde {\mathcal E}_j^{m}\setminus (\tilde {\mathcal E}_j^{m})_{\tilde x_j^m-}^x\big)-\mathcal L^2(E_m\setminus (E_m)_{x_m-}^x)\big)\Big|\leq \frac{1}{m}.
\end{split}
\end{equation}
Let $(j_m)_{m\in\N}$ be an increasing sequence of integer numbers 
such that $j_m\geq J(m)$ for any $m\in\N$.
We choose the diagonal sequence
$\tilde E_m=\tilde{\mathcal E}_{j_m}^m$, $m\in\N$
and prove that 
there exists $\ell>0$
such that
\begin{equation}
\label{bddpr}
\tilde E_m\subset [-\ell,\ell]\times[-\ell, \ell],\quad\text{for any}\quad m\in\N.
\end{equation}
Let $\tilde x_m=\tilde x_{j_m}^m$. We have
\begin{equation}
\label{boundperc}
\sup\{\P((\tilde E_m)_{\tilde x_m-}^x) : m\in\N\}<\infty,\quad\text{and}
\quad
\sup\{\P(\tilde E_m\setminus (\tilde E_m)_{\tilde x_m-}^x) : m\in\N\}<\infty,
\end{equation}
since
$\P((\tilde E_m)_{\tilde x_m-}^x)
+\P(\tilde E_m\setminus(\tilde E_m)_{\tilde x_m-}^x) \leq 2C_{MP}+4h_1+4h_2+2$.

We prove that the sequence $\tilde x_m$ is bounded. 
Let $\tilde \phi_m$ be the profile function of $\tilde E_m$
and assume by contradiction that
$x_m\to\infty$ as $m\to\infty$. In this case, 
by the representation formula \eqref{eq:REPR}, we have
\[
\P((\tilde E_m)_{\tilde x_m-}^x)
=\int_{0}^{\tilde x_m}\sqrt{\tilde \phi_m(x)^2+|x|^{2\alpha}}\;dx
\geq \int_0^{\tilde{x}_m}|x|^\alpha=\frac{\tilde x_m^{\alpha+1}}{\alpha+1}\to\infty,
\quad m\to\infty
\]
which contradicts \eqref{boundperc}.
In the same way, we can see that, if $r_m$ is 
such that $\tilde E_m\subset(\tilde E_m)_{r_m-}^x$,
the sequence $(r_m)_{m\in\N}$ is bounded.

Now, we show that there exists $L\geq 0$ such that
$\tilde  E_m\subset(\tilde  E_m)_{L-}^y$.
Suppose by contradiction that
for any $L\geq 0$, there exists $m=m(L)\in\N$ such that
$(\tilde E_m)_{\tilde x_m-}^x
\setminus (\tilde E_m)_{L-}^y\neq \emptyset$.
By convexity of $(\tilde E_m)_{\tilde x_m-}^x$, 
we can equivalently
assume that for any $L\geq 0$ there exists $j(L)\geq 0$ such that
\begin{equation}\label{limite}
\tilde \phi_m(0)>L \text{ for }m\geq m(L),
\end{equation}
We write for $x\in(0,\tilde x_m)$
\[
\tilde \phi_m(x)=-\int_x^{\tilde x_m}\tilde \phi_m'(\xi)\;d\xi=\int_{x}^{\tilde x_m}|\tilde \phi_m'(\xi)|\;d\xi,
\]
then
\[
\tilde \phi_m(0)=\lim_{x\to0}\int_x^{\tilde x_m}|\tilde \phi_m'(\xi)|\;d\xi=\int_0^{\tilde x_m}|\tilde \phi_m'(\xi)|\;d\xi
\]
which implies, by \eqref{limite}
\[
\lim_{m\to\infty}\int_0^{\tilde x_m}|\tilde \phi_m'(\xi)|\;d\xi=\lim_{m\to\infty}\tilde \phi_m(0)=\infty.
\]
Therefore
\[
\P((\tilde E_m)_{\tilde x_m-}^x)=4\int_0^{\tilde x_m}\sqrt{(\tilde \phi_m'(x))^2+x^{2\alpha}}\;dx\geq\int_0^{\tilde x_m}|\tilde \phi_m'(x)|\;dx\to\infty\text{ as }m\to\infty,
\]
which is in contradiction with 
\eqref{boundperc}.
Similarly, we exclude the case that for any $L>0$
$\big(\tilde E_m\setminus (\tilde E_m)_{\tilde x_m-}^x\big)
\setminus(\tilde E_m)_{L-}^y\neq \emptyset.$

Thanks to \eqref{boundperc}, by
the compactness theorem for $BV_\alpha$ functions
(see \cite[Thm 1.28 ]{GN}), 
there exists
a set $E_\infty$ which is
the $L^1_{loc}$-limit of $\tilde E_m$
as $m\to\infty$.
By \eqref{bddpr}, 
convergence 
$\chi_{\tilde{E}_m}\to\chi_{E_\infty}$ is in $L^1(\R^2)$.
Moreover, since the sequence $(\tilde x_m)_{m\in\N}$ is bounded, 
we let 
$x_\infty\geq 0$ be the limit up to subsequences of $\tilde x_m$
as $m\to\infty$.
We have, by \eqref{feder},
\[
\mathcal L^2((E_\infty)_{\tilde x_\infty-}^x)
=\lim_{m\to\infty}\mathcal L^2((\tilde E_m)_{\tilde x_m-}^x)
=\lim_{m\to\infty}\mathcal L^2(( E_m)_{x_m-}^x)=v_1,
\]
and
\[\mathcal L^2(E_\infty\setminus (E_\infty)_{\tilde x_\infty-}^x)
=\lim_{m\to\infty}\mathcal L^2(\tilde E_m\setminus (\tilde E_m)_{\tilde x_m-}^x)
=\lim_{m\to\infty}\mathcal L^2( E_m\setminus ( E_m)_{x_m-}^x)=v_2.
\]
Now, since $(\tilde E_m)_{\tilde x_m-}^x$,
is convex, we can choose a representative for $E_\infty$ such that
$(E_\infty)_{x_\infty-}^x$,
is convex.
By boundedness of the sequence $\tilde E_m$,
let $y_\infty^\pm\geq 0$ be such that
$\tr_{x_\infty\pm}^xE_\infty=[-y_\infty^\pm,y_\infty^\pm]$.
Then, by \eqref{trac} and claim 4 at Step 4, we have
\[
y_\infty^-\geq \lim_{m\to\infty} \tilde q_m^-\geq \lim_{m\to\infty} q_m^-\geq 
\lim_{m\to\infty} y_m^--\frac{1}{m}\geq h_1,
\]
equivalently $y_\infty^+\geq h_2$. Hence $E_\infty\in\mathcal A$.

By the lower semi-continuity of the $\alpha$-perimeter together with \eqref{feder}, 
we have
$\mathscr P_\alpha(E_\infty)\leq \liminf_{m\to\infty}\mathscr{P}_\alpha(E_m)+\frac{2}{m}\leq C_{MP}$.
Moreover, by Remark \ref{rem:lip} applied to $E=E_\infty$, there exists a set $\tilde E_\infty\in \mathcal A(v_1,v_2,h_1,h_2)$ such that $\mathscr P_\alpha(\tilde E_\infty)\leq \mathscr P_\alpha(E_\infty)$, $\tilde E_\infty\setminus(\tilde E_\infty)_{x_\infty-}^x$  is locally Lipschitz and $(\tilde E_\infty)_{x_\infty-}^x=(E_\infty)_{x_\infty-}^x$.
This concludes the proof.
\endproof

\section{Profile of regular solutions}
\label{S:char}

In this section we describe the solutions to the minimal partition problem \eqref{eq:pmini} found in Theorem \ref{EXISTpart}.
If $E\in\mathcal A(v_1,v_2,h_1,h_2)$ with partitioning point $x_0>0$ is a bounded solution to Problem \eqref{eq:pmini} such that $E_{x_0-}^x$ is convex and $E\setminus E_{x_0-}^x$ is locally Lipschitz,
we say that $E$ is a {\em regular solution of the minimal partition problem}.
In this case, writing $E=\{(x,y)\in\R^2 : |y|<f(|x|)\},$
the profile function $f:[0,\infty)\to[0,\infty)$ is {\em decreasing} in $[0,x_0)$ 
and {\em locally Lipschitz} in $[x_0,\infty)$.
In the following 
we sometimes consider $f$ to be extended to an even function
defined on the whole $\R$, and we still call it the profile function of $E$.

\subsection{Differential equations for the profile function}
\label{SS:ODEpart}


\begin{prop}
\label{prop:eqdiff}
Let $v_1,v_2,h_1,h_2\geq 0$ and $E\subset\R^2$ be a regular solution of the minimal partition problem \eqref{eq:pmini} for $\mP$ in the class $\mathcal A=\mathcal A(v_1,v_2,h_1,h_2)$, with partitioning point $x_0\geq0$.
Then, writing 
\[
E=\{(x,y)\in\R^2 : |y|<f(x)\},
\]
the even function $f$ satisfies
\begin{subequations}
\begin{alignat}{3}
\label{eqdiffcent}
&f'(x)=-\frac{\mathrm{sgn}{x}\; c|x|^{\alpha+1}}{\sqrt{1-cx^2}}\quad 
&\text{if }\quad |x|<x_0,
\\
\label{eqdiffext+}
&f'(x)=\frac{(kx+d)\; x^\alpha}{\sqrt{1-(kx+d)^2}}\quad
&\text{if }\quad x>x_0.\\
\label{eqdiffext-}
&f'(x)=\frac{(kx-d)\; |x|^\alpha}{\sqrt{1-(kx-d)^2}}\quad
&\text{if }\quad x<-x_0
\end{alignat}
\end{subequations}
for some constants $c\geq 0$, $k,\ d\in\R$.
\end{prop}
\proof
By boundedness of the regular minimizer $E$, let
$r_0=\inf\{r>0 : E\subset E_{r-}^x\}<\infty$.

We prove equation \eqref{eqdiffcent}.
For $\psi_1\in C_c^\infty(0,x_0)$
with $\int\psi_1=0$,
and $\varepsilon\in\R$, consider 
the function
$x\mapsto f(|x|)+\varepsilon \psi_1(|x|)$, $x\in\R$,
and define the set
\[
E_\varepsilon =
\{(x,y)\in\R^2 : |y| < f(|x|)+\varepsilon \psi_1(|x|)\}\in\mathcal A.
\]
By the Representation formula for the $\alpha$-perimeter \eqref{eq:REPR}, let
\[
p_1(\varepsilon)=\P((E_\varepsilon)_{x_0-}^x)
=4\Big\{\int_0^{x_0}\sqrt{(f'+\varepsilon\psi')^2+|x|^{2\alpha}}\;dr
+\lim_{x\to x_0^-}f(x)\Big\}.
\]
By minimality of $E$, we then have
\[
\begin{split}
0&=\left.p_1'(\varepsilon)\right|_{\varepsilon=0}
=4\int_{0}^{x_0}\left.\dfrac{d}{d\varepsilon}
\big(\sqrt{(f'+\varepsilon \psi_1')^2+x^{2\alpha}}\big)\right|_{\varepsilon=0}\;dx\\
&=4\int_{0}^{x_0}\dfrac{f'(x)\psi_1'(x)}{\sqrt{f'^2(x)+x^{2\alpha}}}\;dx=
-4\int_{0}^{x_0}\dfrac{d}{dx}\left(\dfrac{f'}{\sqrt{f'^2+x^{2\alpha}}}\right)\;\psi_1(x)\;dx.
\end{split}
\]
By arbitrariness of $\psi_1$, we
deduce the following second order ordinary differential equation, holding for
some $C\in\R$
\begin{equation}
\label{odeprel}
\dfrac{d}{dx}\left(\dfrac{f'(x)}{\sqrt{f'(x)^2+x^{2\alpha}}}\right)=C
\quad \text{for a.e. }
0<x<x_0.
\end{equation}
The normal form of \eqref{odeprel} is
\begin{equation}
\label{normalprel}
f''(x)=\frac{\alpha f'(x)}{x} +\frac{C}{x^{2\alpha}}(f'(x)^2+x^{2\alpha})^{\frac{3}{2}}.
\end{equation}
Now, since $E$ is $x$-symmetric, the function $f$ is even, hence $f'$ is odd
and $f''$ is even. This allows us to extend \eqref{odeprel} to $|x|<x_0$.
Integrating \eqref{odeprel} around 0, we obtain existence of
a constant $d\in\R$ such that for some $\delta>0$,
\[
\dfrac{f'(x)}{\sqrt{f'(x)^2+x^{2\alpha}}}=Cx+d\quad \text{ for }|x|<\delta.
\]
Since $f'$ is odd we deduce that $d=0$, in fact
for $|x|<\delta$
\[
Cx+d=\dfrac{f'(x)}{\sqrt{f'^2(x)-x^{2\alpha}}}=-\dfrac{f'(-x)}{\sqrt{f'^2(-x)+(-x)^2\alpha}}=-(C(-x)+d)=Cx-d.
\] 
Hence \eqref{odeprel} reads
\[
\dfrac{f'(x)}{\sqrt{f'^2(x)+x^{2\alpha}}}=Cx \text{ for }|x|<\delta,
\]
which implies, by monotonicity of $f$, that $C<0$. Letting $c=-C>0$, we hence get the following ordinary differential equation for $f$:
\[
f'(x)=-\mathrm{sgn}(x)\dfrac{c|x|^{\alpha+1}}{\sqrt{1-c^2x^2}}\quad\text{for }|x|<\delta.
\]
A solution to the latter equation can be extended up to 
$(-1/c,1/c)$. This implies
$0<x_0\leq1/c$ and \eqref{eqdiffcent} is proved.

To prove \eqref{eqdiffext+} and \eqref{eqdiffext-}, 
we proceed in the same way, considering a function
$\psi_2\in C_c^\infty(x_0,r_0)$, 
with $\int \psi_2=0$ and the associated perturbation
$f+\eta \psi_2$ for $\eta \in\R$. 
The set $E_\eta=\{(x,y)\in\R^2 : |y|<f(|x|)+\varepsilon \psi_2(|x|)\}$ 
is inside the class $\mathcal A$, hence, as in the previous case, minimality of
$E$ leads to
\[\dfrac{d}{dx}\left(\dfrac{f'(x)}{\sqrt{f'(x)^2+x^{2\alpha}}}\right)=k
\quad \text{for}\quad
x_0<|x|<r_0
\]
for some constant $k\in\R$.
Let $x_0<x<r_0$. An integration between $x_0$ and $x$ shows that,
letting
\[
d=\lim_{x\to x_0^+}\frac{f'(x)}{\sqrt{f'(x)^2+x^{2\alpha}}}-kx_0,
\]
we have
\[
\dfrac{f'(x)}{\sqrt{f'(x)^2+x^{2\alpha}}}=kx+d\quad\text{for }\quad x_0<x<r_0
\]
which is equivalent to \eqref{eqdiffext+}. In particular, $|kx+d|<1$
for $x_0<x<r_0$.

Analogously, for any $x\in(-r_0,-x_0)$, an integration between
$x$ and $-x_0$ shows that
\[
\dfrac{f'(x)}{\sqrt{f'(x)^2+x^{2\alpha}}}=kx-d\quad\text{for }\quad -r_0<x<-x_0,
\]
which leads to \eqref{eqdiffext-} and $|kx-d|<1$ for $-r_0<x<-x_0$.
\endproof

\subsection{Traces of regular solutions}
\label{SS:traces}

In this section, we study traces of minimizers. 
What we expect is that if $E\in\mathcal A(v_1,v_2,h_1,h_2)$ is
a regular solution of the minimal partition problem with partitioning point $x_0$,
then
\[
\tr_{x_0}^{-}E=[-h_1,h_1]\quad\text{and }\quad
\tr_{x_0}^{+}E=[-h_2,h_2].
\]
In Proposition \ref{TRpart} we prove the claim for the left trace, under an additional assumption. 
The case of the right trace is equivalent.

\begin{prop}
\label{TRpart}
Given $v_1,v_2,h_1,h_2\geq 0$, let $E\in\mathcal A(v_1,v_2,h_1,h_2)$ 
be a regular solution of Problem \eqref{eq:pmini} for $\mP$ with partitioning point $x_0>0$, and
let $f:[0,\infty)\to[0,\infty)$ be its profile function.
If 
\[
\lim_{x\to x_0^-}f'(x)>-\infty,
\]
then 
\[\tr_{x_0}^{-}E=[-h_1,h_1].\]
\end{prop}

\proof Assume by contradiction that $f(x_0^-)>h_1$, where
$
f(x_0^-)=\lim_{x\to x_0^-}f(x).
$
We show that in this case, there exists a set $F\in\mathcal A$ such that
$\P(F_{x_0-}^x)<\P(E_{x_0-}^x)$, $\P(F\setminus F_{x_0-}^x)=\P(E\setminus E_{x_0-}^x)$, hence $\mP(F)<\mP(E)$, which
is in contradiction with the minimality of $E$.

\begin{figure}[h!!!]
\vspace{-0.5cm}
\centering
\includegraphics[scale=0.35]{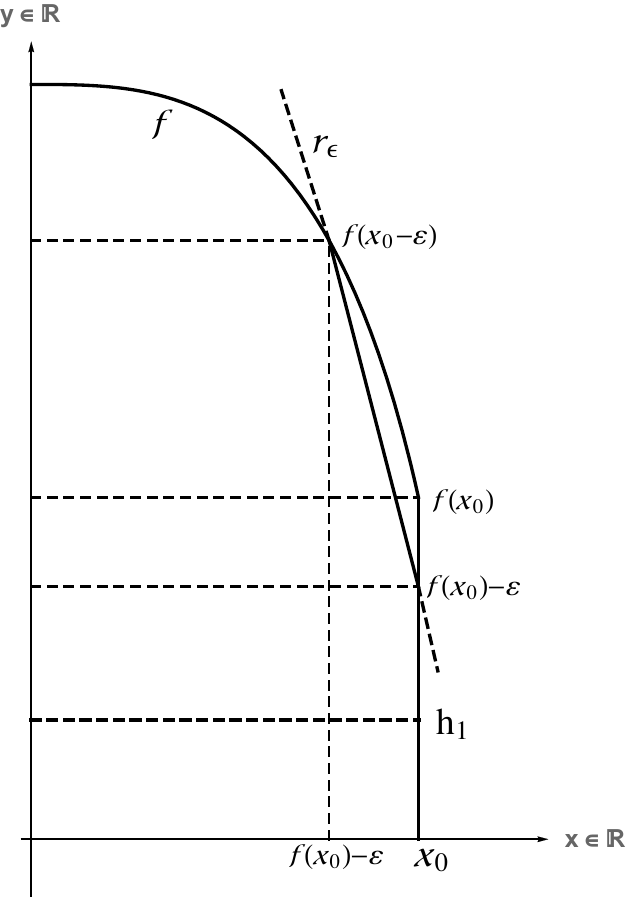}
\vspace{-0.5cm}
\caption{Construction of the set $E_\varepsilon$.}\label{fig:taglio}
\end{figure}
For a small parameter $\varepsilon>0$, let $f_\varepsilon:[0,x_0)\to[0,\infty)$
be the function defined by
\[
f_\varepsilon(x)=\left\{
\begin{array}{l}
f(x),\quad\text{if }0<x<x_0-\varepsilon\\
r_\varepsilon(x),\quad\text{if }x_0-\varepsilon<x<x_0
\end{array}
\right.
\]
where $r_\varepsilon$ is the segment connecting the points
$\big(x_0-\varepsilon,f(x_0-\varepsilon)\big)$ and $(x_0,f(x_0^-)-\varepsilon)$,
i.e.,
\[
r_\varepsilon(x)=m(\varepsilon)(x-x_0)+f(x_0^-)-\varepsilon,
\quad
m(\varepsilon)=\frac{1}{\varepsilon}(f(x_0^-)-\varepsilon-f(x_0-\varepsilon))<0.
\]
We define the set
$E_\varepsilon=\{(x,y)\in\R^2 : |y|<f_\varepsilon(|x|)\}$ and
let $A(\varepsilon)=(\P(E_{x_0-}^x)-\P(E_\varepsilon))/4$,
 $B(\varepsilon)
=(\mathcal L^2(E_{x_0-}^x)-\mathcal L^2(E_\varepsilon))/4$. 
Let $y_\varepsilon=B(\varepsilon)/x_0$.
We claim that for $\varepsilon>0$ small enough the set
\[
F_\varepsilon=(E_\varepsilon+(0,y_\varepsilon))\cup \big([-x_0,x_0]\times[-y_\varepsilon,y_\varepsilon]\big)
\]
obtained by translating $E_\varepsilon$ in the vertical direction of the
quantity $y_\varepsilon$, 
satisfies 
\begin{equation}
\label{claimfe}
\P(F_\varepsilon)< \P(E).
\end{equation}
Since
$\mathcal L^2((F_\varepsilon)_{x_0-}^x)=\mathcal L^2(E_\varepsilon)+4x_0y_\varepsilon=\mathcal L^2(E_{x_0-}^x)$,
this contradicts the minimality of $E$ and the statement follows.

We prove \eqref{claimfe}.
By invariance under vertical translations of the $\alpha$-perimeter, we have
\[
\P(F_\varepsilon)=\P(E_\varepsilon)+4y_\varepsilon=\P(E_{x_0-}^x)-4\Big( A(\varepsilon)-\frac{B(\varepsilon)}{x_0}\Big)
\]
To prove the claim, it is therefore sufficient to show that for $\varepsilon>0$ small
enough
\begin{equation}
\label{claimfin}
x_0A(\varepsilon)> B(\varepsilon).
\end{equation}
This follows by a comparison between the Taylor developments of $A$ and $B$.
We have
\[
\begin{split}
A(\varepsilon)&=\Big\{\int_{x_0-\varepsilon}^{x_0}\sqrt{f'(x)^2+x^{2\alpha}}
-\sqrt{m(\varepsilon)^2+x^{2\alpha}}\;dx +\varepsilon\Big\},
\quad B_\varepsilon=\int_{x_0-\varepsilon}^{x_0}
f(x)-r_\varepsilon(x)\;dx.
\end{split}
\]
Hence 
\begin{equation}
\label{limitr}
\lim_{\varepsilon\to0^+}A(\varepsilon)=\lim_{\varepsilon\to0^+}B(\varepsilon)=0.
\end{equation}
Let $f'(x_0^-)=\lim_{x\to x_0^-}f'(x)<0$
and $f''(x_0^-)=\lim_{x\to x_0^-}f''(x)\leq 0$.
We compute
\[
\lim_{\varepsilon\to0^+}m(\varepsilon)=\lim_{\varepsilon\to0^+}\frac{f(x_0^-)-f(x_0-\varepsilon)}{\varepsilon}-1= f'(x_0^-)-1.
\] 
Now, by the assumption $f'(x_0^-)>-\infty$ we also have
$f''(x_0^-)>-\infty$. Hence 
\[
\begin{split}
m'(0)&=\lim_{\varepsilon\to0^+}m'(\varepsilon)
=\lim_{\varepsilon\to0^+}\frac{1}{\varepsilon}\Big\{f'(x_0-\varepsilon)-\dfrac{1}{\varepsilon}\big[f(x_0^-)-f(x_0-\varepsilon)\big]\Big\}
\\
&=\lim_{\varepsilon\to0^+}\dfrac{1}{\varepsilon}\Big\{f'(x_0^-)-f''(x_0^-)\varepsilon-\dfrac{1}{\varepsilon}\Big[f(x_0^-)-\Big(f(x_0^-)-f'(x_0^-)\varepsilon+\dfrac{f''(x_0^-)}{2}\varepsilon^2\Big)\Big]+o(\varepsilon)\Big\}\\
&=-\frac{f''(x_0^-)}{2}
\end{split}
\]
is positive and finite.
Therefore, 
\[
\begin{split}
A'(\varepsilon)&=1+\sqrt{f'^2(x_0-\varepsilon)+(x_0-\varepsilon)^{2\alpha}}-\sqrt{m(\varepsilon)^2+(x_0-\varepsilon)^{2\alpha}}-
\int_{x_0-\varepsilon}^{x_0}\dfrac{m(\varepsilon)m'(\varepsilon)}{\sqrt{m(\varepsilon^2)+x^{2\alpha}}}\;dx\\
&\geq\sqrt{f'^2(x_0-\varepsilon)+(x_0-\varepsilon)^{2\alpha}}-\sqrt{m(\varepsilon)^2+(x_0-\varepsilon)^{2\alpha}}
\end{split}
\]
that gives
\begin{equation}\label{estdA0}
A'(0)=\lim_{\varepsilon\to0}A'(\varepsilon)
\geq1+\sqrt{f'^2(x_0^-)^2+x_0^{2\alpha}}-\sqrt{(f'(x_0^-)-1)^2+x_0^{2\alpha}}> 0.
\end{equation}
The last inequality is justified by the following: for any  $c<a<0$, and $b\in\R$,
\[
\sqrt{a_2+b_2}-\sqrt{c_2+b_2}> c-a.
\]
We conclude observing that $B'(0)=\lim_{\varepsilon\to0^+}B'(\varepsilon)=0$,
in fact:
\[\begin{split}
B'(\varepsilon)
&=f(x_0-\varepsilon)-r_\varepsilon(x_0-\varepsilon)
+\int_{x_0-\varepsilon}^{x_0}\frac{d}{d\varepsilon}\Big\{f(x)-m(\varepsilon)(x-x_0)-f(x_0^-)+\varepsilon\Big\}\;dx\\
&=f(x_0-\varepsilon)+\varepsilon m(\varepsilon)-f(x_0^-)+\varepsilon
-\int_{x_0-\varepsilon}^{x_0}1-m'(\varepsilon)(x-x_0)\;dx\xrightarrow[\varepsilon\to0]{}0
\end{split}
\]
Then, \eqref{claimfin} follows by \eqref{estdA0} and \eqref{limitr}.
\endproof

\begin{rem}\label{rem:cutno}
Assume that the profile function of a minimizer as in Theorem \ref{EXISTpart}
satisfies 
\[
f'(x_0^-)=\lim_{x\to x_0^-}f'(x)=-\infty.
\]
Using the notation of Proposition \ref{TRpart}, there holds,
for $\varepsilon>0$ small enough:
\[
x_0A(\varepsilon)-B(\varepsilon)=-\frac{x_0^{\alpha+2}}{6\sqrt{2}}\big(\frac{\varepsilon}{x_0}\big)^{\frac{3}{2}}+o(\varepsilon^{3/2})<0\quad\text{ for }\varepsilon<\varepsilon_0.
\]
Hence, the construction proposed in the latter proposition does not apply
to this case.
\end{rem}

\subsection{Center of regular solutions.}
\label{SS:center}
In this section, we show that
we cannot in general expect
a regular
minimizer for the minimal partition problem to be obtained in its central part 
as a dilation of the isoperimetric set $\E^\alpha$. In fact, we show that,
for particular choices of $v_1,v_2,h_1,h_2>0$ and $\alpha\geq 0$, a regular
minimizer $E\in\mathcal A$ with partitioning point $x_0>0$
does not satisfy for any $\lambda>0$
\begin{equation}
\label{concpaul}
E_{x_0-}^x=(\delta_\lambda^\alpha\E^\alpha)_{x_0-}^x.
\end{equation}
\begin{rem} \label{rem:diltra}
Let $E$ be a regular solution of the minimal partition problem with partitioning point $x_0>0$.
Then, by Proposition \ref{prop:eqdiff}, its profile function $f$ is defined on some bounded interval $[0,r_0]$ and it is a locally Lipschitz 
function satisfying in a weak sense the ordinary differential equations \eqref{eqdiffcent}-\eqref{eqdiffext-}. 
By an elementary argument, that is omitted, it follows that 
$f\in C^2([0,x_0))\cap C^2(x_0,r_0)\cap C^2(-r_0,-x_0)$.

Notice that equation \eqref{eqdiffcent} is scale invariant, i.e., given $c_1,c_2\geq 0$ and a solution $g$ to \eqref{eqdiffcent} for $c=c_1$, the function
$ g_\lambda(x)=\lambda^{\alpha+1}g\big(\frac{x}{\lambda}\big)$, for  $\lambda=\frac{c_1}{c_2}$
is a solution to \eqref{eqdiffcent} for $c=c_2$.
In this sense, in \cite[Theorem 3.2]{MM}, the authors show that the unique solution 
to equation \eqref{eqdiffcent} is the function
\begin{equation}
\label{eq:phialpha}
\varphi_\alpha(x)=\int_{\arcsin |x|}^{\frac{\pi}{2}}\sin^{\alpha+1}(t)\;dt,\quad x\in[-1,1],
\end{equation}
obtained integrating \eqref{eqdiffcent} for $c=1$.
In other words, there exist $\lambda>0$ and $y\in\R$ such that
\begin{equation}
\label{diltra}
f(x)=\lambda^{\alpha+1}\varphi_\alpha\big(\frac{x}{\lambda}\big)+y,\quad |x|<x_0.
\end{equation}
\end{rem}
We characterize the parameters $\lambda$ and $y\in\R$ appearing in \eqref{diltra} in terms of the data $h_1=h$, $v_1=v$
in the case of regular solutions of Problem \eqref{eq:pmini}
with $h_2=v_2=0$ and $\alpha=0,1$, see Proposition \ref{ESTh}. 
We deduce that, if $\alpha=1$, 
the translation $y$ is strictly negative (see Remark \ref{rem:noeucl}), hence \eqref{concpaul} does not hold.
On the other hand, in the case when $\alpha=0$, the $\alpha$-perimeter corresponds to the Euclidean perimeter and we prove that regular solutions of the minimal partition problem satisfy \eqref{concpaul}.
\begin{figure}[h]
\subfigure{\includegraphics[scale=0.4]{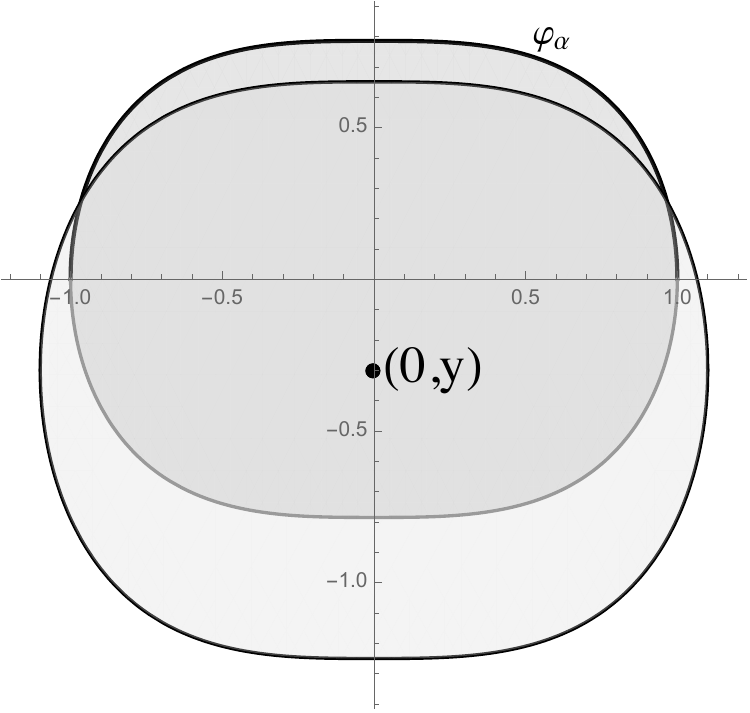}} \hspace{3ex}  
\subfigure{\includegraphics[scale=0.4, trim = 0cm -2.7cm 0cm 0cm]{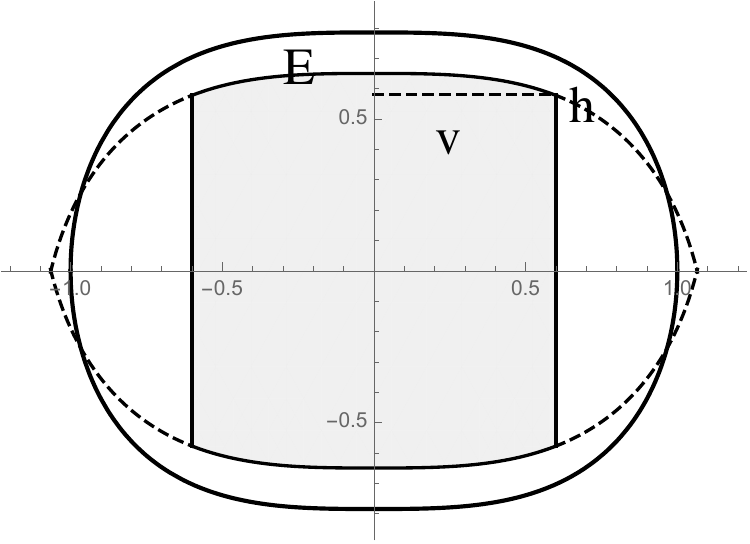}}  
\caption{For $\alpha>0$, a regular solution for the minimal partition problem is not obtained as a dilation of the isoperimetric set $\E^\alpha$ in its central part. Its profile function is in fact the profile of an isoperimetric set vertically translated of a negative quantity $y$.}
\end{figure}

\begin{prop}
\label{ESTh}
Given $\alpha\in\{0,1\}$, $h\geq 0$, and $v>0$,
let  $E\in\mathcal A(v,0,h,0)$ be a regular solution of Problem \eqref{eq:pmini} for $\mP$
with $v_1=v,h_1=h,v_2=h_2=0$,
and partitioning point $x_0>0$, satisfying $\tr_{x_0}^-E=[-h,h]$.
Let $f:[0,x_0]\to[0,\infty)$ be its profile function.
Then there exists $d\in[-1,1]$ such that $f$ solves \eqref{diltra} 
with
\begin{equation}
\label{h(H,V)}
\lambda=\frac{x_0}{d},\quad y=h\Big\{1-\frac{\varphi_\alpha(d)}{d^\alpha\sqrt{1-d^2}}\Big\}.
\end{equation}

\end{prop}
\proof
By Remark \ref{rem:diltra}, there exist $\lambda>0$, $y\in\R$ such that
\[
f(x)=\lambda^{\alpha+1}\varphi_\alpha\big(\frac{x}{\lambda}\big)+y,\quad |x|<x_0.
\]
In particular $\lambda=1/c$.
Let $\beta=f(0)>0$, and 
define
\[
p(\beta,c,x_0)=
\P(E)=\int_0^{x_0}\sqrt{f'^2(t)+t^{2\alpha}}\;dt=\int_0^{x_0}\frac{t^\alpha}{\sqrt{1-(ct)^2}}\;dt
=\frac{1}{c^{\alpha+1}}\int_0^{\arcsin{cx_0}}\!\!\!\!\!\!\!\!\!\!\!\!\!\!\sin^\alpha{\vartheta}\;d\vartheta.
\]
Notice that $p$ can be thought of as a functional depending on $\beta, c, x_0$ where the triple $(\beta, c, x_0)$ identifies a unique minimizer $E\in\A$ as in the statement. In particular, $p$ is independent of $\beta$, as well as $\P$ is independent 
of vertical translations. 
Let $d=c x_0$, that leads to $\lambda=x_0/d$. With a slight abuse of notation, we write $p$ in terms of $d$ and $x_0$ as
\begin{equation}
\label{perv}
p(d,x_0)=x_0^{\alpha+1}g_\alpha(d),\text{ with }g_\alpha(d)=\frac{1}{d^{\alpha+1}}\int_0^{\arcsin{d}}\sin^{\alpha}\vartheta\;d\vartheta.
\end{equation}
We write the volume and trace constraints satisfied by the minimizer $E$ in terms of the parameters $d$ and $x_0$. 
For any $t\in(0,x_0)$, there holds
\begin{equation}
\label{eq:ff}
\begin{split}
f(t)&=\beta+\int_0^tf'(s)\;ds
=\beta-\frac{x_0^{\alpha+1}}{d^{\alpha+1}}\int_0^{\arcsin{(\frac{d}{x_0}t)}}\sin^{\alpha+1}\vartheta\;d\vartheta.
\end{split}
\end{equation}
Hence the trace constraint $f(x_0)=h$ is equivalent to
\begin{equation}
\label{tracecons}
\begin{split}
&\beta=\beta(d,x_0)=h+x_0^{\alpha+1}\sigma_\alpha(d),\text{ with }
\sigma_\alpha(d)=\frac{1}{d^{\alpha+1}}\int_0^{\arcsin{d}}\!\!\!\!\!\!\sin^{\alpha+1}\vartheta\;d\vartheta>0\text{ for }d\in(0,1).
\end{split}
\end{equation}
Plugging \eqref{tracecons} into \eqref{eq:ff} we get
\begin{equation}
\label{f(lambda)}
\begin{split}
y&=f(\lambda)=h+x_0^{\alpha+1} b_\alpha(d),\text{ with }
b_\alpha(d)
=-\frac{\varphi_\alpha(d)}{d^{\alpha+1}}
<0\text{ for }d\in(0,1).\end{split}
\end{equation}
Moreover, the volume constraint $\int f=v$ reads
\[\begin{split}
v&
=\Big(h+x_0^{\alpha+1}\sigma_\alpha(d)\Big) x_0-\frac{x_0^{\alpha+2}}{d^{\alpha+2}}\int_0^d\int_0^{\arcsin{t}}\!\!\!\!\!\!\sin^{\alpha+1}\vartheta\;d\vartheta\;dt,
\end{split}
\]
and we write $v=h x_0+x_0^{\alpha+2}G_\alpha(d),$ with 
\begin{equation}
\label{volume}
\begin{split}
G_\alpha(d)&=\frac{1}{d^{\alpha+2}}\Big(d\int_0^{\arcsin{d}}\sin^{\alpha+1}\vartheta\;d\vartheta-\int_0^d\int_0^{\arcsin{t}}\sin^{\alpha+1}\vartheta\;d\vartheta\;dt\Big)\\
&=\frac{1}{d^{\alpha+2}}\int_0^d\int_{\arcsin{t}}^{\arcsin{d}}\sin^{\alpha+1}\vartheta\;d\vartheta\;dt >0\text{ for }d\in(0,1)
\end{split}
\end{equation}
The functional
\[
F(d,x_0)=x_0^{\alpha+2}G_\alpha(d)+hx_0-v.
\]
defines implicitly the constraints of the problem.
Existence of a minimizer together with the Lagrange
Multipliers theorem imply that there exists $\mu\in\R$ such that
$\nabla p(d,x_0)=\mu\nabla F(d,x_0)$, namely
\begin{equation}
\label{syslagr}
\left\{
\begin{array}{l}
\partial_d p =\mu \partial_d F\\
\partial_{x_0}p=\mu \partial_{x_0} F\\
F(d,x_0)=0
\end{array}
\right.
\iff
\left\{
\begin{array}{l}
g_\alpha'(d) x_0^{\alpha+1}=\mu x_0^{\alpha+2}G_\alpha'(d)\\
(\alpha+1)x_0^\alpha g_\alpha(d)=\mu\big( (\alpha+2)x_0^{\alpha+1}G_\alpha(d)+h\big)\\
x_0^{\alpha+2}G_\alpha(d)+hx_0-v=0
\end{array}
\right.
\end{equation}
Recalling the definitions of $g_\alpha$ and $G_\alpha$
in  \eqref{perv} and \eqref{volume},
we write the expressions for the derivatives
\[\begin{split}
g_\alpha'(d)&=-\frac{\alpha+1}{d^{\alpha+2}}\int_0^{\arcsin{d}}\sin^\alpha\vartheta\;d\vartheta+\frac{1}{d\sqrt{1-d^2}}\\
G_\alpha'(d)&= -\frac{\alpha+2}{d^{\alpha+3}}\int_0^d\int_{\arcsin{t}}^{\arcsin{d}}\sin^{\alpha+1}\vartheta\;d\vartheta
+\frac{1}{\sqrt{1-d^2}}\end{split}.\]
A straightforward computation implies that for $\alpha=0,1$ we have
\begin{equation}
\label{vpsi}
g_\alpha'(d)=dG_\alpha'(d),
\end{equation}
or, equivalently:
\begin{equation}
\label{vpsi2}
\begin{split}
0&=\frac{1}{d^{\alpha+2}}\Big\{(\alpha+2)\int_0^d\int_{\arcsin{t}}^{\arcsin{d}}\!\!\!\!\!\!\sin^{\alpha+1}\vartheta\;d\vartheta\;dt
-(\alpha+1)\int_0^{\arcsin{d}}\!\!\!\!\!\!\sin^{\alpha}\vartheta\;d\vartheta\Big\}+\frac{1-d^2}{d\sqrt{1-d^2}}\\
&=(\alpha+2)G_\alpha(d)-\frac{1}{d}(\alpha+1)g_\alpha(d)+\frac{1}{d}\sqrt{1-d^2}.
\end{split}
\end{equation}

From \eqref{vpsi}, the first equation in system \eqref{syslagr}
gives
\[
\mu=\frac{g_\alpha'(d)}{x_0G_\alpha'(d)}=\frac{d}{x_0}=\frac{1}{\lambda}.
\]
Plugging $\mu$ into the second equation of \eqref{syslagr}
we obtain
\[
(\alpha+1)x_0^{\alpha+1}g_\alpha(d)=d\{(\alpha+2)x_0^{\alpha+1}G_\alpha(d)+h\},
\]
hence, using \eqref{vpsi2},
\begin{equation}
\label{x_0^a+1}
\begin{split}
x_0^{\alpha+1}&=\frac{dh}{(\alpha+1)g_\alpha(d)-d(\alpha+2)G_\alpha(d)}=\frac{dh}{\sqrt{1-d^2}}.
\end{split}
\end{equation}

We are left to compute $y=f(\lambda)=f(x_0/d)$ with $x_0=x_0(h,d)$
given by \eqref{x_0^a+1}.
Expression \eqref{f(lambda)} for $y$, combined with
\eqref{x_0^a+1} gives
\[\begin{split}
y&=f\Big(\frac{x_0}{d}\Big)=h+x_0^{\alpha+1}b_\alpha(d)
=h-\frac{dh}{\sqrt{1-d^2}}\frac{\varphi_\alpha(d)}{d^{\alpha+1}}\\
&=-\frac{h}{d^{\alpha}\sqrt{1-d^2}}\Big\{\varphi_\alpha(d)-d^{\alpha}\sqrt{1-d^2}\Big\}
\end{split}\]
which concludes the proof.
\endproof

\begin{rem}\label{rem:noeucl} 
Let $\alpha\geq 0$ and $E=\{(x,y)\in\R^2 : |y|<f(|x|)\}$ be a regular minimizer of Problem \eqref{eq:pmini}
with $f$ as in \eqref{diltra}. We deduce by \eqref{h(H,V)} that $y=0$ if and only if $\alpha=0$. 
In fact, for any $\alpha\geq 0$
the function $d\mapsto\varphi_\alpha(d)-d^\alpha\sqrt{1-d^2}$
is 0 at $d=1$ and, if $\alpha>0$, it is strictly monotone decreasing since
\[
(\varphi_\alpha(d)-d^\alpha\sqrt{1-d^2})'=-\frac{d^{\alpha+1}}{\sqrt{1-d^2}}-\alpha d^{\alpha-1}\sqrt{1-d^2}+\frac{d^{\alpha+1}}{\sqrt{1-d^2}}<0.
\]
Hence if $\alpha>0$, $\varphi_\alpha(d)-d^\alpha\sqrt{1-d^2}>0$ for $0<d<1$. In particular, $y<0$.
On the other hand, if $\alpha=0$, $\varphi_\alpha(d)=\sqrt{1-d^2}$, that leads to $y=0$.

This implies
that the central part of Euclidean solutions of Problem \eqref{eq:pmini} are 
portions of isoperimetric sets lying in some stripe $\{|x|<x_0\}$, 
while this property fails to hold in the Grushin plane with $\alpha=1$.

\end{rem}

\appendix

\section{Traces of Schwarz symmetric sets}
\label{app:traces}

For a set $E\in\S^*_y$ and a point $x_0\in\R$, the notion of {\em trace of $E$ at $x_0$} can be defined thanks to the following Lemma.
\begin{lem}
\label{lem:traces}
Let $E\in\S^*_y$ and let $x_0\in\R$.
Then there exist $y^+,y^-\geq0$ such that if $T^+=[-y^+,y^+]$  and $T^-=[-y^-,y^-]$,
there holds
\[
\lim_{x\to x_0^\pm}\int_\R|\chi_E(x,y)-\chi_{T^\pm}(y)|\;dy=0.
\]
\end{lem}
\proof
We prove the statement for the limit as $x\to x_0^-$.
Let $u\in C^1(\R^2)$ and $x_1,x_2\in(-\infty,x_0)$. Consider the $\alpha$-gradient
of $u$, $D_\alpha u=(\partial_x u, |x|^\alpha\partial_y u)$. We have
\[
\begin{split}
\int_\R\big(u(x_2,y)-u(x_1,y)\big)dy
&=\int_\R\int_{x_1}^{x_2}\!\!\partial_x u(\xi,y)\;d\xi dy\leq\int_{(x_1,x_2)\times\R}\!\!\!\!\!\!\!\!\!\!\!\!\!\!\!\!|\partial_x u|(\xi,y)\;d\xi dy\leq |D_\alpha u|((x_1,x_2)\times\R).
\end{split}
\]
By the approximation theorem for $BV_\alpha$-functions, see \cite{FSSC},
the last inequality can be extended to $u\in BV_\alpha(\R^2)$
and for $u=\chi_E$ we get
\[
\int_\R\big(\chi_E(x_2,y)-\chi_E(x_1,y)\big)dy\leq\P(E;(x_1,x_2)\times\R),
\]
It hence follows that for every $\varepsilon>0$ there exists $\delta>0$,
such that
\[\|\chi_E(x_2,\cdot)-\chi_E(x_1,\cdot)\|_{L^1(\R)}\leq\varepsilon\text{ for } x_0-\delta<x_1<x_2<x_0.
\]
We deduce existence of 
a function $u\in L^1(\R)$ which is the limit of $\chi_E(x,\cdot)$
as ${x\to x_0^-}$.
Moreover, since for any $x\in\R$, the section $E_x=\{y\in\R : (x,y)\in E\}$ is a
real interval centered at zero, then
$u=\chi_{T^-}$, for a symmetric interval $T^-=[-y^-,y^-]$ for some $y^->0$.
\endproof
\begin{defn}[Traces of Schwarz symmetric sets]
\label{def:traces}
Let $E\in\S^*_y$ be a set  with 
finite $\alpha$-perimeter
and let $x_0\in\R$. 
The interval $T^-$ (resp. $T^+$) defined in Lemma \ref{lem:traces}
is called the {\em left (resp. right) trace of $E$ at $x_0$}
and it is denoted by $\tr_{x_0}^{-}$ (resp. $\tr_{x_0}^{+}$).
If 
\[
\mathrm{tr}_{x_0-}^x E=\mathrm{tr}_{x_0+}^x E=[-y_0,y_0],
\]
we set $\mathrm{tr}_{x_0}^xE=[-y_0,y_0]$ and we call it the {\em trace of $E$ at $x_0$ in the $x$-direction}.
In this case we say that the set {\em $E$ has trace at $x_0$ in the $x$-direction}.
\end{defn}

\begin{rem}
\label{tracesm}
If $E\in\S_x\cap\S^*_y$ has profile function $f:[0,\infty)\to[0,\infty)$
then left and right traces at $x_0>0$ can be computed as:
\[\tr_{x_0\pm}^xE=[-y_0^\pm, y_0^\pm]
\quad\text{with}\quad
\lim_{x\to x_0^\pm}f(x)=y_0^\pm.
\]
In fact, by definition of left and right traces, we have
\[
\begin{split}
0&=\!\!\lim_{x\to x_0^\pm}\!\!\int_\R\!\!|\chi_{E}(x,y)-\chi_{[-y_0^\pm,y_0^\pm]}(y)|\;dy\!
=\!\!\!\lim_{x\to x_0^\pm}\mathcal L^1((E)_x\triangle[-y_0^\pm,y_0^\pm])
=\!2\!\!\lim_{x\to x_0^\pm}|f(x)- y_0^\pm|.
\end{split}
\]
\end{rem}


\section*{Acknowledgments}
\small
The author would like to thank Prof. Roberto Monti from the University of Padova for suggesting 
her the problem and for the precious help.
The author acknowledges the support of the G.N.A.M.P.A. projects
{\em Problemi isoperimetrici e teoria geometrica della misura in spazi metrici},
and {\em Variational Problems and Geometric Measure Theory in Metric Spaces}.

\end{document}